\documentclass[11pt]{article}
\usepackage[top=1.25in, bottom=1.25in, left=1.25in, right=1.25in]{geometry}
\usepackage{amssymb}
\usepackage{latexsym}
\usepackage{amsmath}
\usepackage{amsthm}
\usepackage{mathtools}
\usepackage{tikz-cd}
\usepackage{tikz}
\usepackage{mathrsfs}
\usepackage{stmaryrd}
\usepackage{enumitem} 
\usepackage{hyperref}
\usepackage{graphicx}
\usepackage{caption}
\usepackage{multicol}
\usepackage{subcaption}
\usepackage[mathscr]{eucal}
\usepackage[title]{appendix}
\usepackage{mathrsfs}
\usepackage{bbm}
\usepackage{bm}
\usepackage{booktabs}
\usepackage{color,soul}
\usepackage{stackengine, scalerel}
\usepackage[font=footnotesize]{caption}
\usepackage{circuitikz}
\usetikzlibrary{decorations.markings, positioning, quotes, graphs,graphs.standard, arrows, arrows.meta, angles}

\DeclareMathOperator{\im}{im}

\graphicspath{ {./images/} }
\hypersetup{
    colorlinks,
    citecolor=blue,
    filecolor=red,
    linkcolor=blue,
    urlcolor=blue
}

\theoremstyle{definition}
\newtheorem{definition}{Definition}[section]

\theoremstyle{definition}
\newtheorem{example}{Example}[section]

\newcommand{\real}{\mathbb R}
\newcommand{\cplx}{\mathbb C}
\newcommand{\intg}{\mathbb Z}

\newcommand{\lgl}{\langle}
\newcommand{\rgl}{\rangle}

\begin{document}

\title{Persistent Topological Laplacians  -- a Survey}
\author{Xiaoqi Wei$^{1}\footnote{Corresponding author.  
 E-mail: weixiaoq@msu.edu} $, 
 and Guo-Wei Wei$^{1,2,3}$\footnote{Corresponding author. E-mail: weig@msu.edu} \\% 
$^1$ Department of Mathematics, \\
Michigan State University, MI 48824, USA.\\
$^2$ Department of Electrical and Computer Engineering,\\
Michigan State University, MI 48824, USA. \\
$^3$ Department of Biochemistry and Molecular Biology,\\
Michigan State University, MI 48824, USA. \\
}
\date{}

\maketitle

\begin{abstract}

Persistent topological Laplacians constitute a new class of tools in topological data analysis (TDA). 
They are motivated by the necessity to address challenges encountered in persistent homology when handling complex data. 
These Laplacians combines multiscale analysis with topological techniques to characterize the  topological and geometrical features of functions and data. 
Their kernels fully retrieve the topological invariants of corresponding persistent homology, while their non-harmonic spectra provide supplementary information.
Persistent topological Laplacians have demonstrated superior performance over persistent homology in analyzing large-scale protein engineering datasets.
In this survey, we offer a pedagogical review of persistent topological Laplacians formulated in various mathematical settings, including simplicial complexes, path complexes, flag complexes, digraphs, hypergraphs, hyperdigraphs, cellular sheaves, as well as $N$-chain complexes.

\end{abstract}


\newpage
\tableofcontents
\newpage

\section{Introduction}

In recent years, there has been exponential growth in research on topological data analysis (TDA) in data science. 
TDA provides a set of mathematical and computational techniques to extract insightful information from complex, high-dimensional datasets.
The primary goal of TDA is to uncover and understand the underlying spatial features of data, 
that might be difficult to capture with traditional methods from statistics, physics, and other mathematical disciplines. 
The main tools of TDA are adapted from homology theory.
In homology theory, we associate algebraic objects, such as groups, rings, or modules, with geometrical objects such that topological features can be inferred from algebraic objects.
The most basic geometrical objects in homology theory are simplicial complexes, consisting of simplices that can be used to model a variety of interactions in complex data.
When the input is a point cloud, traditional (simplicial) homology theory only captures trivial topological information. 
Hence, it is impossible to deduce the shape of a point cloud solely through the calculation of simplicial homology.
A major breakthrough in overcoming this limitation is the invention of persistent homology \cite{edelsbrunner2008persistent, zomorodian2005computing}. The basic idea is that, instead of focusing solely on the point cloud itself, one constructs a multiscale family of simplicial complexes from the input point cloud, called a filtration, and examines the evolution of these simplicial complexes and their associated homology groups across scales.
The output of persistent homology are arrays of topological invariants computed over various scales, often visualized or represented by persistence diagrams, persistence barcodes \cite{ghrist2008barcodes}, persistence images \cite{adams2017persistence,xia2015multidimensional}, or persistence landscapes \cite{bubenik2015statistical}. 
Persistent homology has proven to be the most important technique in TDA.

Persistent homology has been applied to a wide variety of disciplines, including image processing \cite{bae2017beyond, clough2020topological}, neuroscience \cite{dabaghian2012topological}, computational chemistry \cite{townsend2020representation}, computational biology \cite{cang2017topologynet, cang2018integration, gameiro2015topological, kovacev2016using}, nanomaterials \cite{lee2017quantifying, xia2015persistent}, crystalline materials \cite{jiang2021topological}, and complex networks \cite{horak2009persistent}, among others.
One of the most remarkable applications of persistent homology is the dominant success of topological deep learning (TDL) models in the D3R Grand Challenges, a global competition series in computer-aided drug design \cite{nguyen2020review, nguyen2019mathematical}. Another notable achievement is the TDL-facilitated discovery of the evolutionary mechanism of SARS-CoV-2 \cite{chen2020mutations}.
The term topological deep learning was coined in 2017 \cite{cang2017topologynet} and has since become a trending topic in data science and machine learning \cite{papamarkou2024position}. The success of TDL, along with other topology-based machine learning algorithms, has made topological data analysis (TDA) a prominent subject in applied mathematics and data science.

However, since homology theory can only characterize topological spaces up to homotopy equivalence, persistent homology has many limitations when dealing with complex data. 
It neglects certain aspects of shape evolution in a filtration that might be important in applications. For example, in a filtration, the zeroth Betti number stops changing once all points are connected, even though the connectivity of simplicial complexes is evolving.
Additionally, for a heterogeneous 1-dimensional cycle, persistent homology does not account for the cycle's composition or the number of points in the cycle. 
These features are particularly important for complex data, such as those encountered in biological sciences. 
Persistent Laplacians \cite{lieutier2014persistent, wang2020persistent} address some of these challenges by introducing a multiscale version of combinatorial Laplacians.
Roughly speaking, Laplacians are matrices whose spectra encapsulate both topological and non-topological information. 
A persistent Laplacian is defined on a pair of simplicial complexes, and its kernel is isomorphic to the corresponding persistent homology group. 
This means that the harmonic spectra of persistent Laplacians can fully recover the barcodes output by persistent homology, and non-harmonic spectra captures additional information about the input filtration.
%Therefore, given a multiscale family of simplicial complexes, persistent Laplacians capture additional details about shape evolution through their non-harmonic spectra.

In spectral graph theory, the graph Laplacian or Kirchhoff matrix is extensively studied \cite{chung1997spectral}. 
Given a graph, the number of zero eigenvalues of its graph Laplacian is equal to the number of connected components of the graph \cite{horak2013spectra}.  
Besides the number of connected components, many properties of a graph is related to the graph Laplacian, such as the relation of Fiedler value and the graph connectivity \cite{butler2010small}. 
However, graph  Laplacian  only accounts for pairwise interactions. Indeed, the graph Laplacian can be seen as a special case, i.e., the first one,  of a series of  combinatorial Laplacians introduced by Eckmann in 1944 \cite{eckmann1944harmonische}, which are defined for each dimension on a simplicial complex. 
It is well-known that the kernel of a combinatorial Laplacian is isomorphic to the corresponding simplicial homology group \cite{eckmann1944harmonische}.  

The relationship between Laplacians and homology has long been discovered in many different contexts.  
On a differentiable manifold, the de Rham-Hodge theory states that the kernel of a Hodge Laplacian is isomorphic to the corresponding de Rham cohomology group.  The discretization of Hodge Laplacians can be achieved by the discrete exterior calculus \cite{desbrun2006discrete,dodziuk1976finite} and finite element exterior calculus \cite{arnold2006finite}. The associated Helmholtz-Hodge decomposition has wide-spread applications in various fields \cite{bhatia2012helmholtz}. 
The Hodge Laplacians on graphs was discussed in \cite{lim2020hodge}. The similarity and difference between combinatorial Laplacians on  simplicial complexes and  Hodge Laplacians on differentiable manifolds was carefully examined in the literature  \cite{ribando2024combinatorial}. 
These Laplacians have found a wide variety  of applications in science and engineering, including   
ranking \cite{jiang2011statistical,ma20111,xu2012hodgerank},
graphics and imaging \cite{hirani2003discrete, ma20111, tong2003discrete},
  games and traffic flows \cite{candogan2011flows},
  deep learning \cite{bronstein2017geometric},
  data representations \cite{chui2018representation},
dimension reduction \cite{perea2018multiscale},
  denoising \cite{schaub2018flow},
 object synchronization \cite{gao2021geometry},
 link prediction \cite{benson2018simplicial},
  sensor network coverage \cite{zhang2018distributed},
	generalizing effective resistance to simplicial complexes \cite{kook2018simplicial},
 cryo-electron microscopy \cite{ye2017cohomology},
 brain networks \cite{lee2019harmonic},
and    biological interactions  \cite{schaub2020random}.
 
A multiscale formulation of Hodge Laplacians on manifolds was introduced in 2019 \cite{chen2019evolutionary}. The resulting evolutionary de Rham-Hodge theory  can be viewed as  persistent Hodge Laplacians.  Both persistent Hodge Laplacians and persistent (combinatorial)  Laplacians  are persistent topological Laplacians (PTLs) that extend the scope and capability of TDA.  In the most general sense, any method that utilizes multiscale topological Laplacians to quantitatively characterize the  topological/geometrical shapes  of point cloud data or differentiable manifolds can be thought of as a persistent topological Laplacian approach.

Persistent Laplacians have been studied extensively in the past few years \cite{gulen2023generalization, liu2023algebraic, memoli2022persistent}. 
In addition to differential manifolds and simplicial complexes, 
persistent topological Laplacians have also been formulated on many other mathematical settings, such as flag complexes \cite{jones2023persistent}, digraphs \cite{wang2023persistent}, cellular sheaves \cite{wei2021persistent}, hypergraphs \cite{liu2021persistent}, and hyperdigraphs \cite{chen2023persistent}.
Computational algorithms  \cite{dong2024faster, memoli2022persistent}, including a software package \cite{wang2021hermes},  have been developed for computing persistent Laplacians. 
Persistent Laplacian approaches  have been applied to protein-ligand binding prediction \cite{meng2021persistent}, interactomic network modeling \cite{du2023multiscale},  gene expression analysis \cite{cottrell2023plpca}, deep mutational scanning \cite{chen2023topological},  and SARS-CoV-2 variant analysis \cite{wei2023persistent}.  The advantage of persistent Laplacians over persistent homology was demonstrated with a collection of 34 datasets in protein engineering \cite{qiu2023persistent}.   The power of persistent Laplacians has been exemplified by its successful forecasting of the emerging dominant SARS-CoV-2 variants \cite{chen2022persistent}.

Both persistent homology and the aforementioned persistent topological Laplacians are constructed based on the properties of chain complexes. 
It is possible to define Mayer homology for the more general $N$-chain complexes \cite{mayer1942new}, where $N$ is an integer. 
Recently, Shen et al. have introduced persistent Mayer homology and persistent Mayer Laplacian  \cite{shen2023persistent} to further extend persistent homology and the aforementioned persistent topological Laplacians to $N$-chain complexes, offering a new development in TDA.

While there are numerous reviews and monographs about persistent homology \cite{carlsson2009topology,edelsbrunner2008persistent,otter2017roadmap,pun2022persistent}, there is no review on persistent topological Laplacians. The primary goal of this survey is to introduce the notion of persistent topological Laplacians to a wider audience and to facilitate the further development on the subject.   
In this survey, we will first introduce the basics of persistent homology, 
and then discuss the theory of persistent Laplacians and some of its recent advances. 
The presentation of mathematics in this paper is  pedagogical, and we hope the survey is more accessible to researchers from diverse backgrounds.

\section{Mathematical preliminaries}

\subsection{Simplicial complexes and  homology}

\begin{figure}[htbp]
    \centering
    \begin{subfigure}{0.30\textwidth}
        \centering
        \begin{tikzpicture}[decoration={markings, mark= at position 0.75 with {\arrow{triangle 45}[scale=4]}}, every edge quotes/.style = {auto, font=\footnotesize}]      
            \coordinate (0) at (0,0);
            \coordinate (1) at (2, 0);
            \coordinate (2) at (1, 1.732);
            \fill[black!20!white] (0) -- (1) -- (2) -- cycle;
            %\draw[->] (1,0.866) arc (90:360:0.288);% syntax (starting point coordinates) arc (starting angle:ending angle:radius)
            \draw (0) -- (1) -- (2) -- (0);
            \node at (-0.5,0) {$0$};
            \node at (2.5,0) {$1$};
            \node at (1,2.232) {$2$};
            \draw[fill=black] (0) circle (2pt);
            \draw[fill=black] (1) circle (2pt);
            \draw[fill=black] (2) circle (2pt);
        \end{tikzpicture}
        \caption{}
        \label{simpleGraph1a}
    \end{subfigure}
    \begin{subfigure}{0.30\textwidth}
        \centering
        \begin{tikzpicture}[decoration={markings, mark= at position 0.75 with {\arrow{triangle 45}[scale=4]}}]      
            \coordinate (0) at (0,0);
            \coordinate (1) at (2, 0);
            \coordinate (2) at (1, 1.732);
            \draw [postaction={decorate}] (0) -- (1);
            \draw [postaction={decorate}] (1) -- (2);
            \draw [postaction={decorate}] (0) -- (2); 
            \node at (-0.5,0) {$0$};
            \node at (2.5,0) {$1$};
            \node at (1,2.232) {$2$};
            \draw[fill=black] (0) circle (2pt);
            \draw[fill=black] (1) circle (2pt);
            \draw[fill=black] (2) circle (2pt);
        \end{tikzpicture}
        \caption{}
        \label{triangle}
    \end{subfigure}
    \begin{subfigure}{0.30\textwidth}
        \centering
        \begin{tikzpicture}[decoration={markings, mark= at position 0.75 with {\arrow{triangle 45}[scale=4]}}]      
            \coordinate (0) at (0,0);
            \coordinate (1) at (2, 0);
            \coordinate (2) at (1, 1.732);
            \draw [postaction={decorate}] (0) -- (1);
            \draw [postaction={decorate}] (1) -- (2);
            \node at (-0.5,0) {$0$};
            \node at (2.5,0) {$1$};
            \node at (1,2.232) {$2$};
            \draw[fill=black] (0) circle (2pt);
            \draw[fill=black] (1) circle (2pt);
            \draw[fill=black] (2) circle (2pt);
        \end{tikzpicture}
        \caption{}
        \label{01_12}
    \end{subfigure}
    \caption{(a) The simplicial complex $\{\{0\}, \{1\}, \{2\}, \{0,1\}, \{0,2\}, \{1,2\}, \{0,1,2\}\}$. (b) The simplicial complex $\{0,1,2, 01,02,12\}$. Arrows emphasize that vertices are ordered.
    (c) The simplicial complex $\{0,1,2,01,12\}$. }
\end{figure}
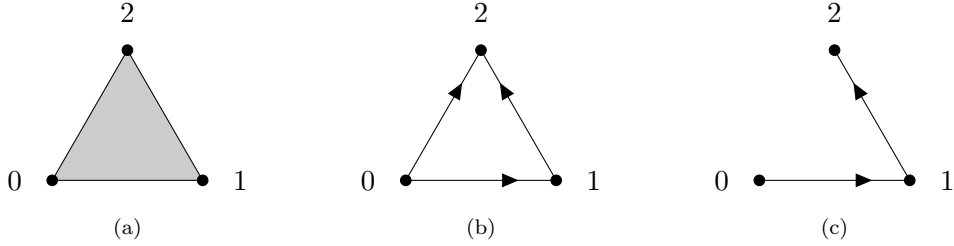

%In discrete mathematics, a graph only describes pairwise interactions between nodes. To describe higher order interactions, we can employ simplicial complexes.  
Given a finite set $V$, a \emph{simplicial complex} $X$ is a collection of subsets of $V$, 
such that if a set $\sigma$ is in $X$, then any subset of $\sigma$ is also in $X$. 
A set $\sigma$ that consists of $q+1$ elements is referred to as a \emph{$q$-simplex}. 
If $\sigma$ is a subset of $\tau$, then we say that $\sigma$ is a face of $\tau$ and denote it by $\sigma \leqslant \tau$.
The definition of a simplicial complex may seem abstract, 
but it is closely related to geometry. 
A $q$-simplex can be realized as the convex hull of $q+1$ points in a real coordinate space, 
so it is possible to construct a polyhedron from a simplicial complex if simplices are glued properly.
For example, suppose $X$ is the power set of $\{0,1,2\}$, 
we can identify $X$ with a triangle whose vertices are labeled as $\{0,1,2\}$ (Figure \ref{simpleGraph1a}).
On the other hand, many geometrical objects can be sliced properly so as to give rise to a simplicial complex.
We always designate a fixed ordering of vertices in a simplicial complex\footnote{The choice of ordering will not affect the resulting homology groups \cite{Munkres1984}.}, 
and require that vertices of any simplices should be ordered according to the fixed ordering.
For example, suppose we use the natural ordering $0<1<2$ for the simplicial complex $\{\{0\},\{1\},\{2\},\{0,1\},\{0,2\},\{1,2\}\}$ (Figure \ref{triangle}), 
then we must not write the simplex $\{0,1\}$ as $\{1,0\}$.
To emphasize that a simplex $\{v_0, \dots, v_q\}$ is ordered, we will use notation $[v_0, \dots, v_q]$ or $v_0 \dots v_q$.

We now introduce a more abstract definition.
A simplicial complex $X$ gives rise to a sequence of vector spaces and linear maps, collectively referred to as a \emph{simplicial chain complex}
\[
    \begin{tikzcd}[column sep = large]
        \centering
    \cdots \arrow[r, "\partial_{3}^X"] & 
    C_{2}(X) \arrow[r, "\partial_{2}^X"] & 
    C_{1}(X)  \arrow[r, "\partial_1^X"] & 
    C_{0}(X) \arrow[r] & 
    0.
    \end{tikzcd} 
\]
The chain group $C_q(X)$ is the real vector space generated by $q$-simplices and the boundary operator $\partial_q$ is a linear map such that 
\begin{align*}
    \partial_q [v_{a_0}, \dots, v_{a_q}] = \sum_i (-1)^i[v_{a_0}, \dots, \hat{v}_{a_i}, \dots, v_{a_q}].
\end{align*}
where the symbol $\hat{v}_{a_i}$ means that $\hat{v}_{a_i}$ is deleted. 
An element of $C_q(X)$ is called a $q$-chain, and by definition it is a linear combination of $q$-simplices. 
Sometimes it is intuitive to regard a $q$-chain as a function mapping a $q$-simplex to its coefficient. 
The coefficients $(-1)^i$ ensure that $\partial_{q}\partial_{q+1}=0$, so the $q$-th homology group $H_q=\ker \partial_{q}/\im \partial_{q+1}$ is well-defined.
The dimension of the homology group $H_q$ is referred to as the \emph{$q$-th Betti number}, 
which is often described as counting the number of $q$-dimensional ``holes'' in a simplicial complex.
It is not always clear what a high-dimensional hole represents in a simplicial complex; nevertheless, the main idea is that homology groups extract quantitative topological information about a simplicial complex.
\begin{example}
    The simplicial complex $X = \{0,1,2,01,02,12\}$ (Figure \ref{triangle}) has only two chain groups $C_0$ and $C_1$, and one boundary map $\partial_1$, represented by the matrix
    \begin{align*}
        \bordermatrix{
            ~  & 01  & 12 & 02  \cr
            0 & -1 & 0 & -1 \cr
            1 & 1  & -1 & 0 \cr
            2 & 0 & 1 & 1 \cr
        }
    \end{align*} 
    if we identify any real-valued function $f_1: \{01,12,02\} \to \real$ with the column vector 
    $(f_1(01), f_1(12), f_1(02))^T$, 
    and any real valued function $f_0: \{0, 1, 2\} \to \real$ with the column vector $(f_0(0), f_0(1), f_0(2))^T$. 
    We can see that $\partial_1 f_1 = f_0$ if and only if $f_0(0)=-f_1(01)-f_1(02), f_0(1)=f_1(01)-f_1(12)$, and $f_0(2)=f_1(12)+f_1(02)$. 
    Since $C_2 = 0$, the homology group $H_1$ is $\ker \partial_1$ and $f_1 \in H_1$ implies $f_1(01)=-f_1(02)=f_1(12)$. 
    
    For the simplicial complex $Y = \{0,1,2,01,12\}$ (Figure \ref{01_12}), the matrix representation of $\partial_1$ is 
    \begin{align*}
    \bordermatrix{
        ~  & 01  & 12   \cr
        0 & -1 & 0  \cr
        1 & 1  & -1 \cr
        2 & 0 & 1  \cr
    }
    \end{align*} 
    and we can verify that the only $f_1$ that satisfies $\partial_1f_1=0$ is the zero function. 
    The intuition behind the difference between $H_1(X)$ and $H_1(Y)$ is that, in $X$ the edges $\{01,12,02\}$ constitute a close path, 
    while in $Y$ there are no close paths.
    %The reader can try to calculate the homology groups of a cycle graph with $n$ vertices.
    %Indeed, for a simple graph, the dimension of $\ker \partial_1$ is equal to the number of $1$-dimensional holes.
\end{example}

A simplicial chain complex is an example of a \emph{chain complex}. 
The reader only need to konw that a chain complex $(V, d)$ is a sequence of vector spaces and linear morphisms
\[
    \begin{tikzcd}[column sep = large]
        \centering
    \cdots \arrow[r, "d_{3}"] & 
    V_{2} \arrow[r, "d_{2}"] & 
    V_{1} \arrow[r, "d_1"] & 
    V_{0} \arrow[r] & 
    0
    \end{tikzcd} 
\]where $d_qd_{q+1} = 0$. We often assume that each $V_q$ is a finite-dimensional inner product space.

\subsection{Combinatorial Laplacians}

Many simplicial complexes share the same Betti numbers.
In this case we can resort to a class of finer descriptors called \emph{combinatorial Laplacians} to distinguish among different simplicial complexes.
Before we define combinatorial Laplacians, we first need to equip a chain group with an inner product.
The canonical approach is to let the set of $q$-simplices be an orthonormal basis for the $q$-th chain group $C_q$.
Now we can discuss the adjoint of the boundary operator $\partial_q$, denoted by $\partial_q^{\ast}$, and
the $q$-th combinatorial Laplacian $\Delta_q$ \cite{eckmann1944harmonische} is defined by 
\begin{align*}
    \partial_{q+1}\partial_{q+1}^{\ast} + \partial_q^{\ast}\partial_q.
\end{align*}
When $q=0$, since $\partial_0 = 0$, the $0$-th combinatorial Laplacian is just $\partial_{1}\partial_{1}^{\ast}$.
The $q$-th combinatorial Laplacian is a positive semi-definite symmetric operator and only has non-negative eigenvalues.
One fact of linear algebra is that, if $U, V,$ and $W$ are inner product spaces and $f:U \to V$, $g: V\to W$ are two linear morphisms such that $g f=0$,
then $\ker (g^{\ast} g + f f^{\ast}) \cong \ker g/\im f$. 
Therefore, the kernel of the $q$-th combinatorial Laplacian $\Delta_q$ is isomorphic to the $q$-th homology group $H_q$ \cite{eckmann1944harmonische}. 
This property guarantees that we can calculate Betti numbers from the spectra of combinatorial Laplacians.
We can further show that $C_q$ admits a \emph{Hodge decomposition} (a detailed exposition can be found in \cite{lim2020hodge})
\begin{align*}
    C_q = \im \partial_q^{\ast} \oplus \ker \Delta_q \oplus \im \partial_{q+1}.
\end{align*}

\begin{example}
    For a simple graph $(V, E)$, let $f_0$ be a function that maps every vertex to a real number.
    If we view the simple graph as a simplicial complex, then $\partial_1^{\ast}$ maps $f_0$ to a real valued function whose domain is $E$.
    The \emph{Dirichlet energy} of $f_0$
    \begin{align*}
        \sum_{v_iv_j\in E}|f_0(v_i)-f_0(v_j)|^2=\lgl \partial_1^{\ast}f_0, \partial_1^{\ast}f_0\rgl = \lgl f_0, \partial_1\partial_1^{\ast} f_0\rgl
    \end{align*}
    measures how $f_0$ varies over $V$. 
    Any $f_0\in \ker \Delta_0 = \ker \partial_1\partial_1^{\ast}$ is a function with zero Dirichlet energy.
    In a connected graph, if $f_0$ has zero Dirichlet energy, then $f_0(a)=f_0(b)$ for any two vertices $a$ and $b$ ($f_0$ is a constant function), 
    because there is always a path that starts from $a$ and ends at $b$. 
    If a graph has more than one connected components, 
    $f_0$ only needs to be constant on any connected components.
    In other words, the dimension of $\ker \Delta_0$ is 
    equal to the number of connected subgraphs. 

    The operator $\Delta_0$ is more commonly known as the graph Laplacian, 
    and there is extensive research studying the relationship between the spectrum of a graph Laplacian and properties of a graph \cite{chung1997spectral}. 
    For a connected graph, it is well known that the minimal nonzero eigenvalue of its graph Laplacian reflects the graph's connectivity \cite{fiedler1973algebraic}.
    Graphs that share the same homology groups may have different graph Laplacians (Figure \ref{HvsL}).
    \begin{figure}[htbp]
        \centering
        \begin{subfigure}{0.30\textwidth}
            \centering
            \begin{tikzpicture}[decoration={markings, mark= at position 0.75 with {\arrow{triangle 45}[scale=4]}}, every edge quotes/.style = {auto, font=\footnotesize}]      
                \coordinate (0) at (1,0);
                \coordinate (1) at (0.5, 0.866);
                \coordinate (2) at (-0.5,0.866);
                \coordinate (3) at (-1,0);
                \coordinate (4) at (-0.5,-0.866);
                \coordinate (5) at (0.5,-0.866);
                \draw (0) -- (1) -- (2) -- (3) -- (4) -- (5);
                \draw[fill=black] (0) circle (2pt);
                \draw[fill=black] (1) circle (2pt);
                \draw[fill=black] (2) circle (2pt);
                \draw[fill=black] (3) circle (2pt);
                \draw[fill=black] (4) circle (2pt);
                \draw[fill=black] (5) circle (2pt);
            \end{tikzpicture}
            \caption{}
            \label{a line}
        \end{subfigure}
        \begin{subfigure}{0.30\textwidth}
            \centering
            \begin{tikzpicture}[decoration={markings, mark= at position 0.75 with {\arrow{triangle 45}[scale=4]}}, every edge quotes/.style = {auto, font=\footnotesize}]      
                \coordinate (0) at (1,0);
                \coordinate (1) at (0.5, 0.866);
                \coordinate (2) at (-0.5,0.866);
                \coordinate (3) at (-1,0);
                \coordinate (4) at (-0.5,-0.866);
                \coordinate (5) at (0.5,-0.866);
                \draw (0) -- (1) -- (2) -- (3) -- (4) -- (5) -- (0);
                \draw[fill=black] (0) circle (2pt);
                \draw[fill=black] (1) circle (2pt);
                \draw[fill=black] (2) circle (2pt);
                \draw[fill=black] (3) circle (2pt);
                \draw[fill=black] (4) circle (2pt);
                \draw[fill=black] (5) circle (2pt);
            \end{tikzpicture}
            \caption{}
            \label{C6}
        \end{subfigure}
        \begin{subfigure}{0.30\textwidth}
            \centering
            \begin{tikzpicture}[decoration={markings, mark= at position 0.75 with {\arrow{triangle 45}[scale=4]}}, every edge quotes/.style = {auto, font=\footnotesize}]      
                \coordinate (0) at (1,0);
                \coordinate (1) at (0,1);
                \coordinate (2) at (-1,0);
                \coordinate (3) at (0,-1);
                \draw (0) -- (1) -- (2) -- (3) -- (0);
                \draw[fill=black] (0) circle (2pt);
                \draw[fill=black] (1) circle (2pt);
                \draw[fill=black] (2) circle (2pt);
                \draw[fill=black] (3) circle (2pt);
            \end{tikzpicture}
            \caption{}
            \label{C4}
        \end{subfigure}
        \caption{Homology can distinguish (a) from (b) and (c), but cannot distinguish between (b) and (c). 
        Laplacians can distinguish among all of them. For a cycle graph with $n$ vertices, 
        the spectrum of the graph Laplacian is $\{2-2\cos (2k\pi/n)\mid k = 1,\dots ,n \}$.}
        \label{HvsL}
    \end{figure}
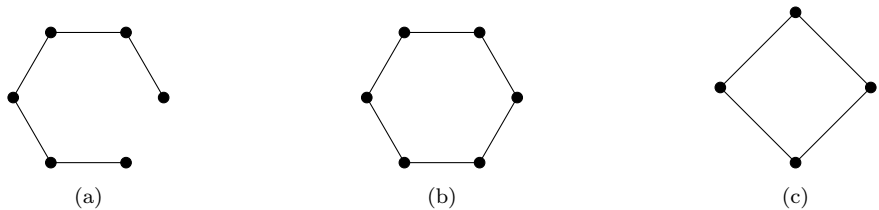
\end{example}

\subsection{Filtration and persistent homology}

So far we have introduced the elementary theory of simplicial complexes, but we have not explained how it is related to point cloud data.
A point cloud is a set of finitely many points $P = \{v_0, \dots, v_n\}$ in a Euclidean space. 
Usually the geometrical structure of a point cloud is related to some non-geometrical properties of the object the point cloud represents,
and a good understanding of the ``geometry'' of the point cloud is important. 
Here the first problem is what we mean by the word ``geometry''. 
A naive notion of ``geometry'' is the pairwise distances between each pair of points,
and we can represent such information by a \emph{filtration}, i.e., a nested sequence of simplicial complexes.
One commonly used filtration is the \emph{Vietoris-Rips filtration}:
given a point cloud $P = \{v_0, \dots, v_n\}$ and a parameter $d\in \real$, 
$X_d$ is a simplicial complex such that the simplex $\{v_{a_0},\dots,v_{a_q}\} \in X_d$ 
if and only if the euclidean distance between $v_{a_i}$ and $v_{a_j}$ is at most $d$ for any $0\leq i < j \leq q$.
By varying $d$, we obtain finitely many distinct simplicial complexes, each of which characterizes the shape of the point cloud at a different scale.
Homology groups or combinatorial Laplacians of each $X_d$ will change as $d$ varies, providing a characterization of the point cloud.
\begin{example}
    Now we build a Vietoris-Rips filtration from the point cloud $\{x=(1,0), y=(0,1), z=(-1,0), w=(0,-1)\} \subset \real^2$ shown in \ref{fourPoints}.
    %For simplicity we only construct $1$-simplices between points.
    When $d=0$, there are no edges in $X_d$. When $d=\sqrt{2}$, $X_d$ changes for the first time and becomes $\{x,y,z,w,xy,yz,zw,xw\}$. 
    If $d$ goes from $\sqrt{2}$ to $2$, $X_2 = X_{\sqrt{2}} \cup \{xz,yw,yzw,xzw,xyw,xyz,xyzw\}$. As $d$ becomes bigger, $X_d$ contains more and more simplices.
    Let's examine $H_1$ for $d=0,\sqrt{2},2$. $H_1(X_0) = 0$ as there is no close path, and $H_1(X_{\sqrt{2}}) = 1$ because of four newly born edges.
    When $d=2$, $H_1 = 0$ since the close path is filled by high-dimensional simplices. 
    \begin{figure}[htbp]
        \centering
        \begin{subfigure}{0.25\textwidth}
            \centering
            \begin{tikzpicture}[decoration={markings, mark= at position 0.75 with {\arrow{triangle 45}[scale=4]}}, every edge quotes/.style = {auto, font=\footnotesize}]      
                \coordinate (0) at (1,0);
                \coordinate (1) at (0,1);
                \coordinate (2) at (-1,0);
                \coordinate (3) at (0,-1);
                %\draw (0) -- (1) -- (2) -- (3) -- (0);
                \draw[fill=black] (0) circle (2pt);
                \draw[fill=black] (1) circle (2pt);
                \draw[fill=black] (2) circle (2pt);
                \draw[fill=black] (3) circle (2pt);
            \end{tikzpicture}
            \caption{}
            \label{fourPoints}
        \end{subfigure}
        \begin{subfigure}{0.25\textwidth}
            \centering
            \begin{tikzpicture}[decoration={markings, mark= at position 0.75 with {\arrow{triangle 45}[scale=4]}}, every edge quotes/.style = {auto, font=\footnotesize}]      
                \coordinate (0) at (1,0);
                \coordinate (1) at (0,1);
                \coordinate (2) at (-1,0);
                \coordinate (3) at (0,-1);
                \draw (0) -- (1) -- (2) -- (3) -- (0);
                \draw[fill=black] (0) circle (2pt);
                \draw[fill=black] (1) circle (2pt);
                \draw[fill=black] (2) circle (2pt);
                \draw[fill=black] (3) circle (2pt);
            \end{tikzpicture}
            \caption{}
            \label{}
        \end{subfigure}
        \begin{subfigure}{0.25\textwidth}
            \centering
            \begin{tikzpicture}[decoration={markings, mark= at position 0.75 with {\arrow{triangle 45}[scale=4]}}, every edge quotes/.style = {auto, font=\footnotesize}]      
                \coordinate (0) at (1,0);
                \coordinate (1) at (0,1);
                \coordinate (2) at (-1,0);
                \coordinate (3) at (0,-1);
                \fill[black!20!white] (0) -- (1) -- (2) -- (3) -- cycle;
                \draw (0) -- (1) -- (2) -- (3) -- (0);
                \draw (0) -- (2);
                \draw (1) -- (3);
                \draw[fill=black] (0) circle (2pt);
                \draw[fill=black] (1) circle (2pt);
                \draw[fill=black] (2) circle (2pt);
                \draw[fill=black] (3) circle (2pt);
            \end{tikzpicture}
            \caption{}
            \label{}
        \end{subfigure}  
        \caption{(a) $X_0 = \{x,y,z,w\}$. (b) $X_{\sqrt{2}} = \{x,y,z,w,xy,yz,zw,xw\}$. (c) $X_2 = X_{\sqrt{2}} \cup \{xz,yw,yzw,xzw,xyw,xyz,xyzw\}$}
        \label{filtration}
    \end{figure}
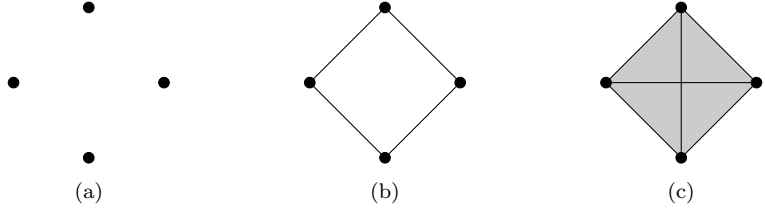 
\end{example}
%It is straightforward to verify that $X_{s} \subset X_{t}$ if $s \leq t$. 

Besides calculating homology groups or combinatorial Laplacians for each $X_t$ in a filtration,
we can also calculate the \emph{persistent homology} to quantify how topological features of a smaller complex $X_{s}$ persist in $X_{t}$. 
Suppose $X$ and $Y$ are two simplicial complexes and $X\subset Y$, then we have the following diagram (dashed arrows indicates inclusion maps $\iota$)
\[
    \begin{tikzcd}[column sep = large]
        C_{q+1}(X) \arrow[r, "\partial_{q+1}^X"] \arrow[d, hook, dashed] 
          & C_{q}(X) \arrow[r, "\partial_q^X"] \arrow[d, hook, dashed] 
            & C_{q-1}(X) \arrow[d, dashed, hook]
            \\
        C_{q+1}(Y) \arrow[r, "\partial_{q+1}^Y"] 
          & C_q(Y) \arrow[r, "\partial_q^Y"] 
            & C_{q-1}(Y). 
    \end{tikzcd}
\]
Since $\im \partial_{q+1}^{Y}$ is larger than $\im \partial_{q+1}^{X}$,
some $q$-dimensional ``holes'' in $X$ might be filled because of $\im \partial_{q+1}^{Y}$.
The $q$-dimensional ``holes'' in $X$ that persists in $Y$ are $\ker \partial_q^X / \im \partial_{q+1}^{Y}$, 
but since 
$\im \partial_{q+1}^{Y}$ is not necessarily a subspace of $\ker \partial_q^X$, 
the proper expression should be 
\begin{align*}
    \ker \partial_q^X / (\im \partial_{q+1}^{Y} \cap \ker \partial_q^X) = \ker \partial_q^X/(\im \partial_{q+1}^{Y} \cap C_q(X)).
\end{align*} 
This quotient space is called the \emph{$q$-th persistent homology group} of the pair $X \subset Y$, the dimension of which is referred to as the \emph{$q$-th persistent Betti number}.

%general definition of persistent homology 
A more formal understanding of persistent homology is helpful. 
We notice that $\partial_q^Y \iota = \iota \partial_q^X$. In plain words, this means that the boundary of a simplex $\sigma\in X$ is unchanged if we view it as a simplex in $Y$.
In general, for two chain complexes $(V, d^V)$ and $(W, d^W)$, the collection of maps $f_q$ such that $f_{q+1} d_q^{V} = d_{q}^W f_q$ for all $q$ is called a chain map.
A chain map $f$ induces a homomorphism $f_{\bullet}:H_q(V) \to H_q(W)$, and sometimes the image $f_{\bullet}(H_q(V))$ is called the persistent homology group of the chain map $f: (V, d^V)\to (W, d^W)$.

%The information of persistent Betti numbers are often represented by a \emph{persistent diagram} or \emph{barcodes} \cite{computationaltopology}.

\section{Persistent (combinatorial) Laplacians}
\subsection{Persistent Laplacians}

We have shown that the kernel of the $q$-th combinatorial Laplacian is isomorphic to the $q$-th homology group. 
This result has been generalized for persistent homology groups.
\[
    \begin{tikzcd}
        C_{q+1}(X) \arrow[rr, "\partial_{q+1}^X"] \arrow[dd, hook, dashed] 
         &
          & C_{q}(X) \arrow[ld, "(\partial_{q+1}^{X,Y})^{\ast}", shift left] \arrow[rr, "\partial_q^X", shift left] \arrow[dd, hook, dashed] 
           &
            & C_{q-1}(X) \arrow[ll, "(\partial_q^X)^{\ast}", shift left] 
            \\
         & C_{q+1}^{X,Y} \arrow[ld, hook, dashed] \arrow[ru, "\partial_{q+1}^{X,Y}", shift left] 
          & 
           &
            &\\
        C_{q+1}(Y) \arrow[rr, "\partial_{q+1}^Y"] 
         & 
          & C_q(Y) 
           &
            &  
    \end{tikzcd}
\]
Given two complexes $X\subset Y$, let $C_{q+1}^{X,Y}$ be the subspace 
\begin{align*}
    \{c \in C_{q+1}(Y) \mid \partial_{q+1}^Y(c) \in C_q(X) \}
\end{align*}
of $C_{q+1}(Y)$ and $\partial_{q+1}^{X,Y}: C_{q+1}^{X,Y}\to C_q(X)$ the restriction of $\partial^Y_{q+1}$ onto $C_{q+1}^{X,Y}$,
then the persistent homology group $\ker \partial_q^X / (\im \partial_{q+1}^{Y} \cap C_q(X))$ is equal to
$\ker \partial_q^X /\im \partial_{q+1}^{X, Y}$. 
Since $C_{q+1}^{X,Y}$ inherits the inner product structure from $C_{q+1}(Y)$,
and $\partial_q^X \partial_{q+1}^{X,Y} = 0$, 
if we define the \emph{$q$-th persistent Laplacian} $\Delta_{q}^{X,Y}: C_q(X)\to C_q(X)$ \cite{lieutier2014persistent, wang2020persistent} by
\begin{align}
    \partial_{q+1}^{X,Y} (\partial_{q+1}^{X,Y})^{\ast} + (\partial^X_{q})^{\ast} \partial^X_{q}
\end{align}
(where $\partial_{q+1}^{X,Y} (\partial_{q+1}^{X,Y})^{\ast}$ is called the up persistent Laplacian, denoted by $\Delta_{q,+}^{X,Y}$, and $(\partial^X_{q})^{\ast} \partial^X_{q}$ the down persistent Laplacian, denoted by $\Delta_{q, -}^{X,Y}$),
we can prove the \emph{persistent Hodge theorem} 
\begin{align*}
    \ker \Delta_{q}^{X,Y} \cong \frac{\ker \partial_q^X}{\im \partial_{q+1}^{Y} \cap C_q(X)} 
\end{align*}
and the persistent Hodge decomposition 
\begin{align*}
    C_q(X) = \im \partial_{q+1}^{X,Y} \oplus \ker\Delta_{q}^{X,Y} \oplus \im (\partial_q^X)^{\ast}
\end{align*}
or equivalently 
\begin{align*}
    C_q(X) = \im \Delta_{q,+}^{X,Y} \oplus \ker\Delta_{q}^{X,Y} \oplus \im \Delta_{q,-}^{X,Y}.
\end{align*}
and the proofs are the same as those for combinatorial Laplacians .
When $X=Y$, $\Delta_q^{X,Y}$ is just a combinatorial Laplacian.
The persistent Hodge theorem implies that information of persistent Betti numbers is included in the spectra of persistent Laplacians,
and additional information can be extracted from nonzero eigenvalues of persistent Laplacians.
Given a point cloud and a filtration $\{X_d, d\in \real\}$ constructed from it, we can calculate persistent Laplacians for 
$X_s\leq X_t$ for a set of preselected $s$ and $t$. 
Employing the information captured in non-harmonic spectra of persistent Laplacians can boost performance of persistent homology-based machine learning models.
In fact even the minimal nonzero eigenvalue of $\Delta^{X_d,X_d}$ already provide a lot of extra information about the point cloud.

\begin{example}
    We illustrate the Vietoris-Rips filtration of a point cloud in Figure \ref{rips}.
    Some results of Laplacian calculation are shown in Figure \ref{res8},
    where $d$ is the diameter (stepsize is 0.02), $\lambda^d_q$ is the minimal nonzero eigenvalue of the $q$-th combinatorial Laplacian of $X_d$,
    and red bars represent homology classes that persist over $d$.   
    The minimal nonzero eigenvalues change at different $d$, indicating the formation of new simplices. 
    \begin{figure}[htbp]
        \centering
        \begin{subfigure}{0.30\textwidth}
            \centering
            \begin{tikzpicture}      
                \coordinate (0) at (0.7,0);
                \coordinate (1) at (1,1);
                \coordinate (2) at (0,0.7);
                \coordinate (3) at (-1,1);
                \coordinate (4) at (-0.7,0);
                \coordinate (5) at (-1,-1);
                \coordinate (6) at (0,-0.7);
                \coordinate (7) at (1,-1);
                \draw[fill=black] (0) circle (2pt);
                \draw[fill=black] (1) circle (2pt);
                \draw[fill=black] (2) circle (2pt);
                \draw[fill=black] (3) circle (2pt);
                \draw[fill=black] (4) circle (2pt);
                \draw[fill=black] (5) circle (2pt);
                \draw[fill=black] (6) circle (2pt);
                \draw[fill=black] (7) circle (2pt);
            \end{tikzpicture}
            \caption{$d=0$}
            \label{}
        \end{subfigure}
        \begin{subfigure}{0.30\textwidth}
            \centering
            \begin{tikzpicture}
                \coordinate (0) at (0.7,0);
                \coordinate (1) at (1,1);
                \coordinate (2) at (0,0.7);
                \coordinate (3) at (-1,1);
                \coordinate (4) at (-0.7,0);
                \coordinate (5) at (-1,-1);
                \coordinate (6) at (0,-0.7);
                \coordinate (7) at (1,-1);
                \draw (0) -- (2) -- (4) -- (6) -- (0);
                \draw[fill=black] (0) circle (2pt);
                \draw[fill=black] (1) circle (2pt);
                \draw[fill=black] (2) circle (2pt);
                \draw[fill=black] (3) circle (2pt);
                \draw[fill=black] (4) circle (2pt);
                \draw[fill=black] (5) circle (2pt);
                \draw[fill=black] (6) circle (2pt);
                \draw[fill=black] (7) circle (2pt);
            \end{tikzpicture}
            \caption{$d=0.99$}
            \label{}
        \end{subfigure}
        \begin{subfigure}{0.30\textwidth}
            \centering
            \begin{tikzpicture}
                \coordinate (0) at (0.7,0);
                \coordinate (1) at (1,1);
                \coordinate (2) at (0,0.7);
                \coordinate (3) at (-1,1);
                \coordinate (4) at (-0.7,0);
                \coordinate (5) at (-1,-1);
                \coordinate (6) at (0,-0.7);
                \coordinate (7) at (1,-1);
                \fill[black!20!white] (0) -- (1) -- (2) -- cycle;
                \fill[black!20!white] (2) -- (3) -- (4) -- cycle;
                \fill[black!20!white] (4) -- (5) -- (6) -- cycle;
                \fill[black!20!white] (0) -- (6) -- (7) -- cycle;
                \draw (0) -- (2) -- (4) -- (6) -- (0);
                \draw (0) -- (1) -- (2) -- (3) -- (4) -- (5) -- (6) -- (7) -- (0);
                \draw[fill=black] (0) circle (2pt);
                \draw[fill=black] (1) circle (2pt);
                \draw[fill=black] (2) circle (2pt);
                \draw[fill=black] (3) circle (2pt);
                \draw[fill=black] (4) circle (2pt);
                \draw[fill=black] (5) circle (2pt);
                \draw[fill=black] (6) circle (2pt);
                \draw[fill=black] (7) circle (2pt);
            \end{tikzpicture}
            \caption{$d = 1.04$}
            %\label{}
        \end{subfigure}
        \begin{subfigure}{0.30\textwidth}
            \centering
            \begin{tikzpicture}
                \coordinate (0) at (0.7,0);
                \coordinate (1) at (1,1);
                \coordinate (2) at (0,0.7);
                \coordinate (3) at (-1,1);
                \coordinate (4) at (-0.7,0);
                \coordinate (5) at (-1,-1);
                \coordinate (6) at (0,-0.7);
                \coordinate (7) at (1,-1);
                \fill[black!20!white] (0) -- (1) -- (2) -- cycle;
                \fill[black!20!white] (2) -- (3) -- (4) -- cycle;
                \fill[black!20!white] (4) -- (5) -- (6) -- cycle;
                \fill[black!20!white] (0) -- (6) -- (7) -- cycle;
                \fill[black!20!white] (0) -- (2) -- (4) -- (6) --  cycle;    
                \draw (0) -- (2) -- (4) -- (6) -- (0);
                \draw (0) -- (1) -- (2) -- (3) -- (4) -- (5) -- (6) -- (7) -- (0);
                \draw (0) -- (4);
                \draw (2) -- (6);
                \draw[fill=black] (0) circle (2pt);
                \draw[fill=black] (1) circle (2pt);
                \draw[fill=black] (2) circle (2pt);
                \draw[fill=black] (3) circle (2pt);
                \draw[fill=black] (4) circle (2pt);
                \draw[fill=black] (5) circle (2pt);
                \draw[fill=black] (6) circle (2pt);
                \draw[fill=black] (7) circle (2pt);
            \end{tikzpicture}
            \caption{$d = 1.40$}
            %\label{}
        \end{subfigure}    
        \begin{subfigure}{0.30\textwidth}
            \centering
            \begin{tikzpicture}
                \coordinate (0) at (0.7,0);
                \coordinate (1) at (1,1);
                \coordinate (2) at (0,0.7);
                \coordinate (3) at (-1,1);
                \coordinate (4) at (-0.7,0);
                \coordinate (5) at (-1,-1);
                \coordinate (6) at (0,-0.7);
                \coordinate (7) at (1,-1);
                \fill[black!20!white] (0) -- (1) -- (2) -- cycle;
                \fill[black!20!white] (2) -- (3) -- (4) -- cycle;
                \fill[black!20!white] (4) -- (5) -- (6) -- cycle;
                \fill[black!20!white] (0) -- (6) -- (7) -- cycle;
                \fill[black!20!white] (0) -- (2) -- (4) -- (6) --  cycle;    
                \draw (0) -- (2) -- (4) -- (6) -- (0);
                \draw (0) -- (1) -- (2) -- (3) -- (4) -- (5) -- (6) -- (7) -- (0);
                \draw (0) -- (4);
                \draw (2) -- (6);
                \draw (0) -- (3);
                \draw (1) -- (4);
                \draw (0) -- (5);
                \draw (2) -- (5);
                \draw (1) -- (6);
                \draw (3) -- (6);
                \draw (2) -- (7);
                \draw (4) -- (7);
                \draw[fill=black] (0) circle (2pt);
                \draw[fill=black] (1) circle (2pt);
                \draw[fill=black] (2) circle (2pt);
                \draw[fill=black] (3) circle (2pt);
                \draw[fill=black] (4) circle (2pt);
                \draw[fill=black] (5) circle (2pt);
                \draw[fill=black] (6) circle (2pt);
                \draw[fill=black] (7) circle (2pt);
            \end{tikzpicture}
            \caption{$d = 1.97$}
            %\label{}
        \end{subfigure}  
        \begin{subfigure}{0.30\textwidth}
            \centering
            \begin{tikzpicture}
                \coordinate (0) at (0.7,0);
                \coordinate (1) at (1,1);
                \coordinate (2) at (0,0.7);
                \coordinate (3) at (-1,1);
                \coordinate (4) at (-0.7,0);
                \coordinate (5) at (-1,-1);
                \coordinate (6) at (0,-0.7);
                \coordinate (7) at (1,-1);
                \fill[black!20!white] (1) -- (3) -- (5) -- (7) -- (1); 
                \draw (0) -- (2) -- (4) -- (6) -- (0);
                \draw (0) -- (1) -- (2) -- (3) -- (4) -- (5) -- (6) -- (7) -- (0);
                \draw (0) -- (4);
                \draw (2) -- (6);
                \draw (0) -- (3);
                \draw (1) -- (4);
                \draw (0) -- (5);
                \draw (2) -- (5);
                \draw (1) -- (6);
                \draw (3) -- (6);
                \draw (2) -- (7);
                \draw (4) -- (7);
                \draw (1) -- (3) -- (5) -- (7) -- (1);
                \draw[fill=black] (0) circle (2pt);
                \draw[fill=black] (1) circle (2pt);
                \draw[fill=black] (2) circle (2pt);
                \draw[fill=black] (3) circle (2pt);
                \draw[fill=black] (4) circle (2pt);
                \draw[fill=black] (5) circle (2pt);
                \draw[fill=black] (6) circle (2pt);
                \draw[fill=black] (7) circle (2pt);
            \end{tikzpicture}
            \caption{$d = 2.00$}
        \end{subfigure}
        \begin{subfigure}{0.30\textwidth}
            \centering
            \begin{tikzpicture}
                \coordinate (0) at (0.7,0);
                \coordinate (1) at (1,1);
                \coordinate (2) at (0,0.7);
                \coordinate (3) at (-1,1);
                \coordinate (4) at (-0.7,0);
                \coordinate (5) at (-1,-1);
                \coordinate (6) at (0,-0.7);
                \coordinate (7) at (1,-1);
                \fill[black!20!white] (1) -- (3) -- (5) -- (7) -- (1); 
                \draw (0) -- (2) -- (4) -- (6) -- (0);
                \draw (0) -- (1) -- (2) -- (3) -- (4) -- (5) -- (6) -- (7) -- (0);
                \draw (0) -- (4);
                \draw (2) -- (6);
                \draw (0) -- (3);
                \draw (1) -- (4);
                \draw (0) -- (5);
                \draw (2) -- (5);
                \draw (1) -- (6);
                \draw (3) -- (6);
                \draw (2) -- (7);
                \draw (4) -- (7);
                \draw (1) -- (5);
                \draw (3) -- (7);
                \draw (1) -- (3) -- (5) -- (7) -- (1);
                \draw[fill=black] (0) circle (2pt);
                \draw[fill=black] (1) circle (2pt);
                \draw[fill=black] (2) circle (2pt);
                \draw[fill=black] (3) circle (2pt);
                \draw[fill=black] (4) circle (2pt);
                \draw[fill=black] (5) circle (2pt);
                \draw[fill=black] (6) circle (2pt);
                \draw[fill=black] (7) circle (2pt);
            \end{tikzpicture}
            \caption{$d = 2.83$}
        \end{subfigure}
        \begin{subfigure}{0.70\textwidth}
            \centering
            \includegraphics[width=0.80\textwidth]{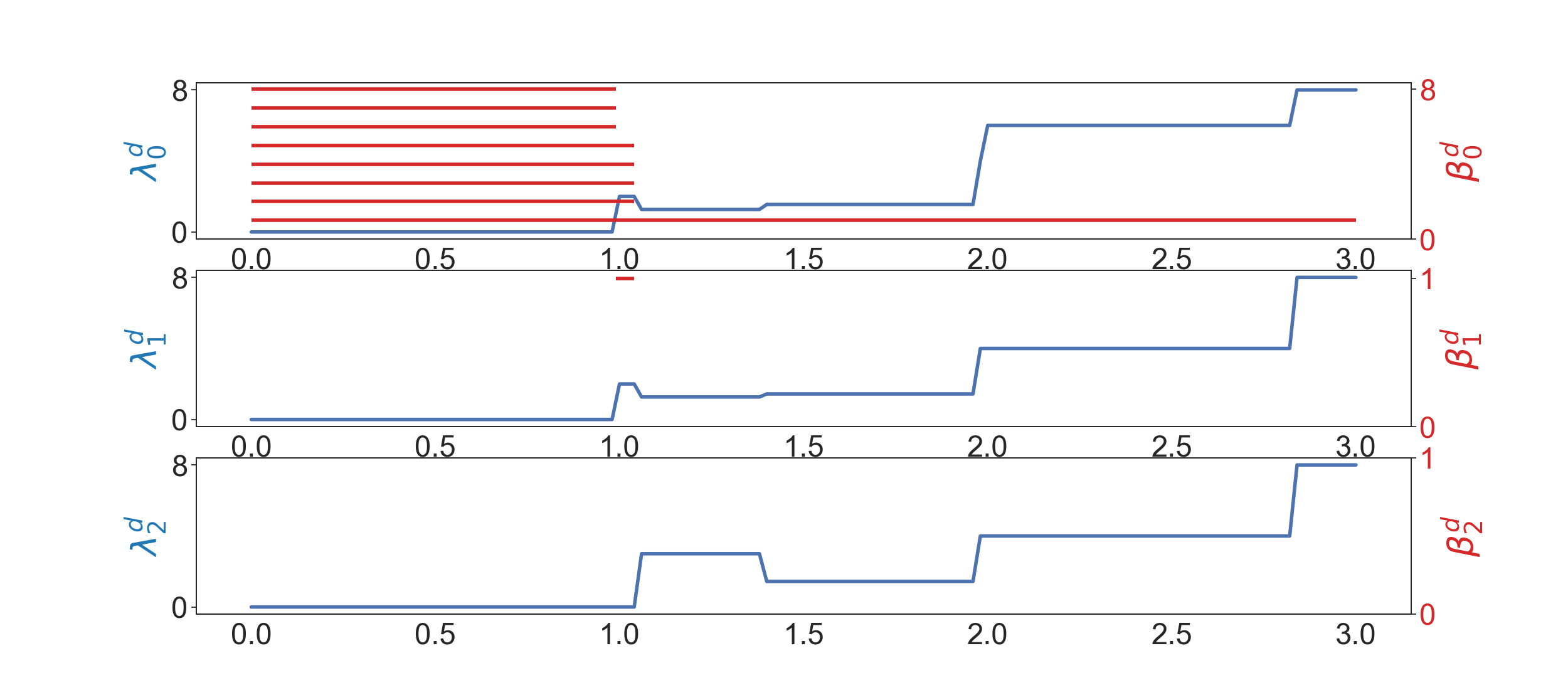}
            %\caption{{\color{red}The results for the Vietoris-Rips filtration constructed from the 8 points.}}
            \caption{}
            \label{res8}
        \end{subfigure}
        \caption{}
        \label{rips}
    \end{figure}
    %Note that $\lambda$ changes more often than $\beta$.
\end{example}

\begin{example}
    \begin{figure}[htbp]
        \centering
        \begin{subfigure}{0.30\textwidth}
            \centering
            \begin{tikzpicture}      
                \coordinate (0) at (1,0);
                \coordinate (1) at (0.866,0.5);
                \coordinate (2) at (0.5, 0.866);
                \coordinate (3) at (0, 1);
                \coordinate (4) at (-0.5, 0.866);
                \coordinate (5) at (-0.866, 0.5);
                \coordinate (6) at (-1, 0);
                \coordinate (7) at (-0.866, -0.5);
                \coordinate (8) at (-0.5, -0.866);
                \coordinate (9) at (0, -1);
                \coordinate (10) at (0.5, -0.866);
                \coordinate (11) at (0.866, -0.5);
                \draw[fill=black] (0) circle (2pt);
                \draw[fill=black] (1) circle (2pt);
                \draw[fill=black] (2) circle (2pt);
                \draw[fill=black] (3) circle (2pt);
                \draw[fill=black] (4) circle (2pt);
                \draw[fill=black] (5) circle (2pt);
                \draw[fill=black] (6) circle (2pt);
                \draw[fill=black] (7) circle (2pt);
                \draw[fill=black] (8) circle (2pt);
                \draw[fill=black] (9) circle (2pt);
                \draw[fill=black] (10) circle (2pt);
                \draw[fill=black] (11) circle (2pt);
            \end{tikzpicture}
            \caption{$d=0$}
            \label{}
        \end{subfigure}
        \begin{subfigure}{0.30\textwidth}
            \centering
            \begin{tikzpicture}
                \coordinate (0) at (1,0);
                \coordinate (1) at (0.866,0.5);
                \coordinate (2) at (0.5, 0.866);
                \coordinate (3) at (0, 1);
                \coordinate (4) at (-0.5, 0.866);
                \coordinate (5) at (-0.866, 0.5);
                \coordinate (6) at (-1, 0);
                \coordinate (7) at (-0.866, -0.5);
                \coordinate (8) at (-0.5, -0.866);
                \coordinate (9) at (0, -1);
                \coordinate (10) at (0.5, -0.866);
                \coordinate (11) at (0.866, -0.5);
                \draw (0) -- (1) -- (2) -- (3) -- (4) -- (5) -- (6) -- (7) -- (8) -- (9) -- (10) -- (11) -- (0);
                \draw[fill=black] (0) circle (2pt);
                \draw[fill=black] (1) circle (2pt);
                \draw[fill=black] (2) circle (2pt);
                \draw[fill=black] (3) circle (2pt);
                \draw[fill=black] (4) circle (2pt);
                \draw[fill=black] (5) circle (2pt);
                \draw[fill=black] (6) circle (2pt);
                \draw[fill=black] (7) circle (2pt);
                \draw[fill=black] (8) circle (2pt);
                \draw[fill=black] (9) circle (2pt);
                \draw[fill=black] (10) circle (2pt);
                \draw[fill=black] (11) circle (2pt);
            \end{tikzpicture}
            \caption{$d=0.52$}
            \label{}
        \end{subfigure}
        \begin{subfigure}{0.30\textwidth}
            \centering
            \begin{tikzpicture}
                \coordinate (0) at (1,0);
                \coordinate (1) at (0.866,0.5);
                \coordinate (2) at (0.5, 0.866);
                \coordinate (3) at (0, 1);
                \coordinate (4) at (-0.5, 0.866);
                \coordinate (5) at (-0.866, 0.5);
                \coordinate (6) at (-1, 0);
                \coordinate (7) at (-0.866, -0.5);
                \coordinate (8) at (-0.5, -0.866);
                \coordinate (9) at (0, -1);
                \coordinate (10) at (0.5, -0.866);
                \coordinate (11) at (0.866, -0.5);
                \fill[black!20!white] (0) -- (1) -- (2) -- cycle;
                \fill[black!20!white] (1) -- (2) -- (3) -- cycle;
                \fill[black!20!white] (2) -- (3) -- (4) -- cycle;
                \fill[black!20!white] (3) -- (4) -- (5) -- cycle;
                \fill[black!20!white] (4) -- (5) -- (6) -- cycle;
                \fill[black!20!white] (5) -- (6) -- (7) -- cycle;
                \fill[black!20!white] (6) -- (7) -- (8) -- cycle;
                \fill[black!20!white] (7) -- (8) -- (9) -- cycle;
                \fill[black!20!white] (8) -- (9) -- (10) -- cycle;
                \fill[black!20!white] (9) -- (10) -- (11) -- cycle;
                \fill[black!20!white] (10) -- (11) -- (0) -- cycle;
                \fill[black!20!white] (11) -- (0) -- (1) -- cycle;
                \draw (0) -- (1) -- (2) -- (3) -- (4) -- (5) -- (6) -- (7) -- (8) -- (9) -- (10) -- (11) -- (0);
                \draw (1) -- (3) -- (5) -- (7) -- (9) -- (11) -- (1);
                \draw (0) -- (2) -- (4) -- (6) -- (8) -- (10) -- (0);
                \draw[fill=black] (0) circle (2pt);
                \draw[fill=black] (1) circle (2pt);
                \draw[fill=black] (2) circle (2pt);
                \draw[fill=black] (3) circle (2pt);
                \draw[fill=black] (4) circle (2pt);
                \draw[fill=black] (5) circle (2pt);
                \draw[fill=black] (6) circle (2pt);
                \draw[fill=black] (7) circle (2pt);
                \draw[fill=black] (8) circle (2pt);
                \draw[fill=black] (9) circle (2pt);
                \draw[fill=black] (10) circle (2pt);
                \draw[fill=black] (11) circle (2pt);
            \end{tikzpicture}
            \caption{$d = 1.00$}
            %\label{}
        \end{subfigure}
        \begin{subfigure}{0.30\textwidth}
            \centering
            \begin{tikzpicture}
                \coordinate (0) at (1,0);
                \coordinate (1) at (0.866,0.5);
                \coordinate (2) at (0.5, 0.866);
                \coordinate (3) at (0, 1);
                \coordinate (4) at (-0.5, 0.866);
                \coordinate (5) at (-0.866, 0.5);
                \coordinate (6) at (-1, 0);
                \coordinate (7) at (-0.866, -0.5);
                \coordinate (8) at (-0.5, -0.866);
                \coordinate (9) at (0, -1);
                \coordinate (10) at (0.5, -0.866);
                \coordinate (11) at (0.866, -0.5);
                \fill[black!20!white] (0) -- (1) -- (2) -- cycle;
                \fill[black!20!white] (1) -- (2) -- (3) -- cycle;
                \fill[black!20!white] (2) -- (3) -- (4) -- cycle;
                \fill[black!20!white] (3) -- (4) -- (5) -- cycle;
                \fill[black!20!white] (4) -- (5) -- (6) -- cycle;
                \fill[black!20!white] (5) -- (6) -- (7) -- cycle;
                \fill[black!20!white] (6) -- (7) -- (8) -- cycle;
                \fill[black!20!white] (7) -- (8) -- (9) -- cycle;
                \fill[black!20!white] (8) -- (9) -- (10) -- cycle;
                \fill[black!20!white] (9) -- (10) -- (11) -- cycle;
                \fill[black!20!white] (10) -- (11) -- (0) -- cycle;
                \fill[black!20!white] (11) -- (0) -- (1) -- cycle;
                \fill[black!20!white] (0) -- (1) -- (2) -- (3) -- cycle;
                \fill[black!20!white] (1) -- (2) -- (3) -- (4) -- cycle;
                \fill[black!20!white] (2) -- (3) -- (4) -- (5) -- cycle;
                \fill[black!20!white] (3) -- (4) -- (5) -- (6) -- cycle;
                \fill[black!20!white] (4) -- (5) -- (6) -- (7) -- cycle;
                \fill[black!20!white] (5) -- (6) -- (7) -- (8) -- cycle;
                \fill[black!20!white] (6) -- (7) -- (8) -- (9) -- cycle;
                \fill[black!20!white] (7) -- (8) -- (9) -- (10) -- cycle;
                \fill[black!20!white] (8) -- (9) -- (10) -- (11) -- cycle;
                \fill[black!20!white] (9) -- (10) -- (11) -- (0) -- cycle;
                \fill[black!20!white] (10) -- (11) -- (0) -- (1) -- cycle;
                \fill[black!20!white] (11) -- (0) -- (1) -- (2) -- cycle;
                \draw (0) -- (1) -- (2) -- (3) -- (4) -- (5) -- (6) -- (7) -- (8) -- (9) -- (10) -- (11) -- (0);
                \draw (1) -- (3) -- (5) -- (7) -- (9) -- (11) -- (1);
                \draw (0) -- (2) -- (4) -- (6) -- (8) -- (10) -- (0);
                \draw (0) -- (3) -- (6) -- (9) -- (0);
                \draw (1) -- (4) -- (7) -- (10) -- (1);
                \draw (2) -- (5) -- (8) -- (11) -- (2);
                \draw[fill=black] (0) circle (2pt);
                \draw[fill=black] (1) circle (2pt);
                \draw[fill=black] (2) circle (2pt);
                \draw[fill=black] (3) circle (2pt);
                \draw[fill=black] (4) circle (2pt);
                \draw[fill=black] (5) circle (2pt);
                \draw[fill=black] (6) circle (2pt);
                \draw[fill=black] (7) circle (2pt);
                \draw[fill=black] (8) circle (2pt);
                \draw[fill=black] (9) circle (2pt);
                \draw[fill=black] (10) circle (2pt);
                \draw[fill=black] (11) circle (2pt);
            \end{tikzpicture}
            \caption{$d = 1.42$}
            %\label{}
        \end{subfigure}    
        \begin{subfigure}{0.30\textwidth}
            \centering
            \begin{tikzpicture}
                \coordinate (0) at (1,0);
                \coordinate (1) at (0.866,0.5);
                \coordinate (2) at (0.5, 0.866);
                \coordinate (3) at (0, 1);
                \coordinate (4) at (-0.5, 0.866);
                \coordinate (5) at (-0.866, 0.5);
                \coordinate (6) at (-1, 0);
                \coordinate (7) at (-0.866, -0.5);
                \coordinate (8) at (-0.5, -0.866);
                \coordinate (9) at (0, -1);
                \coordinate (10) at (0.5, -0.866);
                \coordinate (11) at (0.866, -0.5);
                \fill[black!20!white] (0) -- (1) -- (2) -- cycle;
                \fill[black!20!white] (1) -- (2) -- (3) -- cycle;
                \fill[black!20!white] (2) -- (3) -- (4) -- cycle;
                \fill[black!20!white] (3) -- (4) -- (5) -- cycle;
                \fill[black!20!white] (4) -- (5) -- (6) -- cycle;
                \fill[black!20!white] (5) -- (6) -- (7) -- cycle;
                \fill[black!20!white] (6) -- (7) -- (8) -- cycle;
                \fill[black!20!white] (7) -- (8) -- (9) -- cycle;
                \fill[black!20!white] (8) -- (9) -- (10) -- cycle;
                \fill[black!20!white] (9) -- (10) -- (11) -- cycle;
                \fill[black!20!white] (10) -- (11) -- (0) -- cycle;
                \fill[black!20!white] (11) -- (0) -- (1) -- cycle;
                \fill[black!20!white] (0) -- (1) -- (2) -- (3) -- (4) -- cycle;
                \fill[black!20!white] (1) -- (2) -- (3) -- (4) -- (5) -- cycle;
                \fill[black!20!white] (2) -- (3) -- (4) -- (5) -- (6) -- cycle;
                \fill[black!20!white] (3) -- (4) -- (5) -- (6) -- (7) -- cycle;
                \fill[black!20!white] (4) -- (5) -- (6) -- (7) -- (8) -- cycle;
                \fill[black!20!white] (5) -- (6) -- (7) -- (8) -- (9) -- cycle;
                \fill[black!20!white] (6) -- (7) -- (8) -- (9) -- (10) -- cycle;
                \fill[black!20!white] (7) -- (8) -- (9) -- (10) -- (11) -- cycle;
                \fill[black!20!white] (8) -- (9) -- (10) -- (11) -- (0) -- cycle;
                \fill[black!20!white] (9) -- (10) -- (11) -- (0) -- (1) -- cycle;
                \fill[black!20!white] (10) -- (11) -- (0) -- (1) -- (2) -- cycle;
                \fill[black!20!white] (11) -- (0) -- (1) -- (2) -- (3) -- cycle;
                \draw (0) -- (1) -- (2) -- (3) -- (4) -- (5) -- (6) -- (7) -- (8) -- (9) -- (10) -- (11) -- (0);
                \draw (1) -- (3) -- (5) -- (7) -- (9) -- (11) -- (1);
                \draw (0) -- (2) -- (4) -- (6) -- (8) -- (10) -- (0);
                \draw (0) -- (3) -- (6) -- (9) -- (0);
                \draw (1) -- (4) -- (7) -- (10) -- (1);
                \draw (2) -- (5) -- (8) -- (11) -- (2);
                \draw (0) -- (4) -- (8) -- (0);
                \draw (1) -- (5) -- (9) -- (1);
                \draw (2) -- (6) -- (10) -- (2);
                \draw (3) -- (7) -- (11) -- (3);
                \draw[fill=black] (0) circle (2pt);
                \draw[fill=black] (1) circle (2pt);
                \draw[fill=black] (2) circle (2pt);
                \draw[fill=black] (3) circle (2pt);
                \draw[fill=black] (4) circle (2pt);
                \draw[fill=black] (5) circle (2pt);
                \draw[fill=black] (6) circle (2pt);
                \draw[fill=black] (7) circle (2pt);
                \draw[fill=black] (8) circle (2pt);
                \draw[fill=black] (9) circle (2pt);
                \draw[fill=black] (10) circle (2pt);
                \draw[fill=black] (11) circle (2pt);
            \end{tikzpicture}
            \caption{$d = 1.74$}
            %\label{}
        \end{subfigure}  
        \begin{subfigure}{0.30\textwidth}
            \centering
            \begin{tikzpicture}
                \coordinate (0) at (1,0);
                \coordinate (1) at (0.866,0.5);
                \coordinate (2) at (0.5, 0.866);
                \coordinate (3) at (0, 1);
                \coordinate (4) at (-0.5, 0.866);
                \coordinate (5) at (-0.866, 0.5);
                \coordinate (6) at (-1, 0);
                \coordinate (7) at (-0.866, -0.5);
                \coordinate (8) at (-0.5, -0.866);
                \coordinate (9) at (0, -1);
                \coordinate (10) at (0.5, -0.866);
                \coordinate (11) at (0.866, -0.5);
                \fill[black!20!white] (0) -- (1) -- (2) -- cycle;
                \fill[black!20!white] (1) -- (2) -- (3) -- cycle;
                \fill[black!20!white] (2) -- (3) -- (4) -- cycle;
                \fill[black!20!white] (3) -- (4) -- (5) -- cycle;
                \fill[black!20!white] (4) -- (5) -- (6) -- cycle;
                \fill[black!20!white] (5) -- (6) -- (7) -- cycle;
                \fill[black!20!white] (6) -- (7) -- (8) -- cycle;
                \fill[black!20!white] (7) -- (8) -- (9) -- cycle;
                \fill[black!20!white] (8) -- (9) -- (10) -- cycle;
                \fill[black!20!white] (9) -- (10) -- (11) -- cycle;
                \fill[black!20!white] (10) -- (11) -- (0) -- cycle;
                \fill[black!20!white] (11) -- (0) -- (1) -- cycle;
                \fill[black!20!white] (0) -- (1) -- (2) -- (3) -- (4) -- (5) -- cycle;
                \fill[black!20!white] (1) -- (2) -- (3) -- (4) -- (5) -- (6) -- cycle;
                \fill[black!20!white] (2) -- (3) -- (4) -- (5) -- (6) -- (7) -- cycle;
                \fill[black!20!white] (3) -- (4) -- (5) -- (6) -- (7) -- (8) -- cycle;
                \fill[black!20!white] (4) -- (5) -- (6) -- (7) -- (8) -- (9) -- cycle;
                \fill[black!20!white] (5) -- (6) -- (7) -- (8) -- (9) -- (10) -- cycle;
                \fill[black!20!white] (6) -- (7) -- (8) -- (9) -- (10) -- (11) -- cycle;
                \fill[black!20!white] (7) -- (8) -- (9) -- (10) -- (11) -- (0) -- cycle;
                \fill[black!20!white] (8) -- (9) -- (10) -- (11) -- (0) -- (1) -- cycle;
                \fill[black!20!white] (9) -- (10) -- (11) -- (0) -- (1) -- (2) -- cycle;
                \fill[black!20!white] (10) -- (11) -- (0) -- (1) -- (2) -- (3) -- cycle;
                \fill[black!20!white] (11) -- (0) -- (1) -- (2) -- (3) -- (4) -- cycle;
                \draw (0) -- (1) -- (2) -- (3) -- (4) -- (5) -- (6) -- (7) -- (8) -- (9) -- (10) -- (11) -- (0);
                \draw (1) -- (3) -- (5) -- (7) -- (9) -- (11) -- (1);
                \draw (0) -- (2) -- (4) -- (6) -- (8) -- (10) -- (0);
                \draw (0) -- (3) -- (6) -- (9) -- (0);
                \draw (1) -- (4) -- (7) -- (10) -- (1);
                \draw (2) -- (5) -- (8) -- (11) -- (2);
                \draw (0) -- (4) -- (8) -- (0);
                \draw (1) -- (5) -- (9) -- (1);
                \draw (2) -- (6) -- (10) -- (2);
                \draw (3) -- (7) -- (11) -- (3);
                \draw (0) -- (5) -- (10) -- (3) -- (8) -- (1) -- (6) -- (11) -- (4) -- (9) -- (2) -- (7) -- (0);
                \draw[fill=black] (0) circle (2pt);
                \draw[fill=black] (1) circle (2pt);
                \draw[fill=black] (2) circle (2pt);
                \draw[fill=black] (3) circle (2pt);
                \draw[fill=black] (4) circle (2pt);
                \draw[fill=black] (5) circle (2pt);
                \draw[fill=black] (6) circle (2pt);
                \draw[fill=black] (7) circle (2pt);
                \draw[fill=black] (8) circle (2pt);
                \draw[fill=black] (9) circle (2pt);
                \draw[fill=black] (10) circle (2pt);
                \draw[fill=black] (11) circle (2pt);
            \end{tikzpicture}
            \caption{$d = 1.94$}
        \end{subfigure}
        \begin{subfigure}{0.30\textwidth}
            \centering
            \begin{tikzpicture}
                \coordinate (0) at (1,0);
                \coordinate (1) at (0.866,0.5);
                \coordinate (2) at (0.5, 0.866);
                \coordinate (3) at (0, 1);
                \coordinate (4) at (-0.5, 0.866);
                \coordinate (5) at (-0.866, 0.5);
                \coordinate (6) at (-1, 0);
                \coordinate (7) at (-0.866, -0.5);
                \coordinate (8) at (-0.5, -0.866);
                \coordinate (9) at (0, -1);
                \coordinate (10) at (0.5, -0.866);
                \coordinate (11) at (0.866, -0.5);
                \fill[black!20!white] (0) -- (1) -- (2) -- (3) -- (4) -- (5) -- (6) -- (7) -- (8) -- (9) -- (10) -- (11) -- cycle;
                \draw (0) -- (1) -- (2) -- (3) -- (4) -- (5) -- (6) -- (7) -- (8) -- (9) -- (10) -- (11) -- (0);
                \draw (1) -- (3) -- (5) -- (7) -- (9) -- (11) -- (1);
                \draw (0) -- (2) -- (4) -- (6) -- (8) -- (10) -- (0);
                \draw (0) -- (3) -- (6) -- (9) -- (0);
                \draw (1) -- (4) -- (7) -- (10) -- (1);
                \draw (2) -- (5) -- (8) -- (11) -- (2);
                \draw (0) -- (4) -- (8) -- (0);
                \draw (1) -- (5) -- (9) -- (1);
                \draw (2) -- (6) -- (10) -- (2);
                \draw (3) -- (7) -- (11) -- (3);
                \draw (0) -- (5) -- (10) -- (3) -- (8) -- (1) -- (6) -- (11) -- (4) -- (9) -- (2) -- (7) -- (0);
                \draw (0) -- (6);
                \draw (1) -- (7);
                \draw (2) -- (8);
                \draw (3) -- (9);
                \draw (4) -- (10);
                \draw (5) -- (11);
                \draw[fill=black] (0) circle (2pt);
                \draw[fill=black] (1) circle (2pt);
                \draw[fill=black] (2) circle (2pt);
                \draw[fill=black] (3) circle (2pt);
                \draw[fill=black] (4) circle (2pt);
                \draw[fill=black] (5) circle (2pt);
                \draw[fill=black] (6) circle (2pt);
                \draw[fill=black] (7) circle (2pt);
                \draw[fill=black] (8) circle (2pt);
                \draw[fill=black] (9) circle (2pt);
                \draw[fill=black] (10) circle (2pt);
                \draw[fill=black] (11) circle (2pt);
            \end{tikzpicture}
            \caption{$d = 2.00$}
        \end{subfigure}
        \begin{subfigure}{0.70\textwidth}
            \centering
            \includegraphics[width=0.80\textwidth]{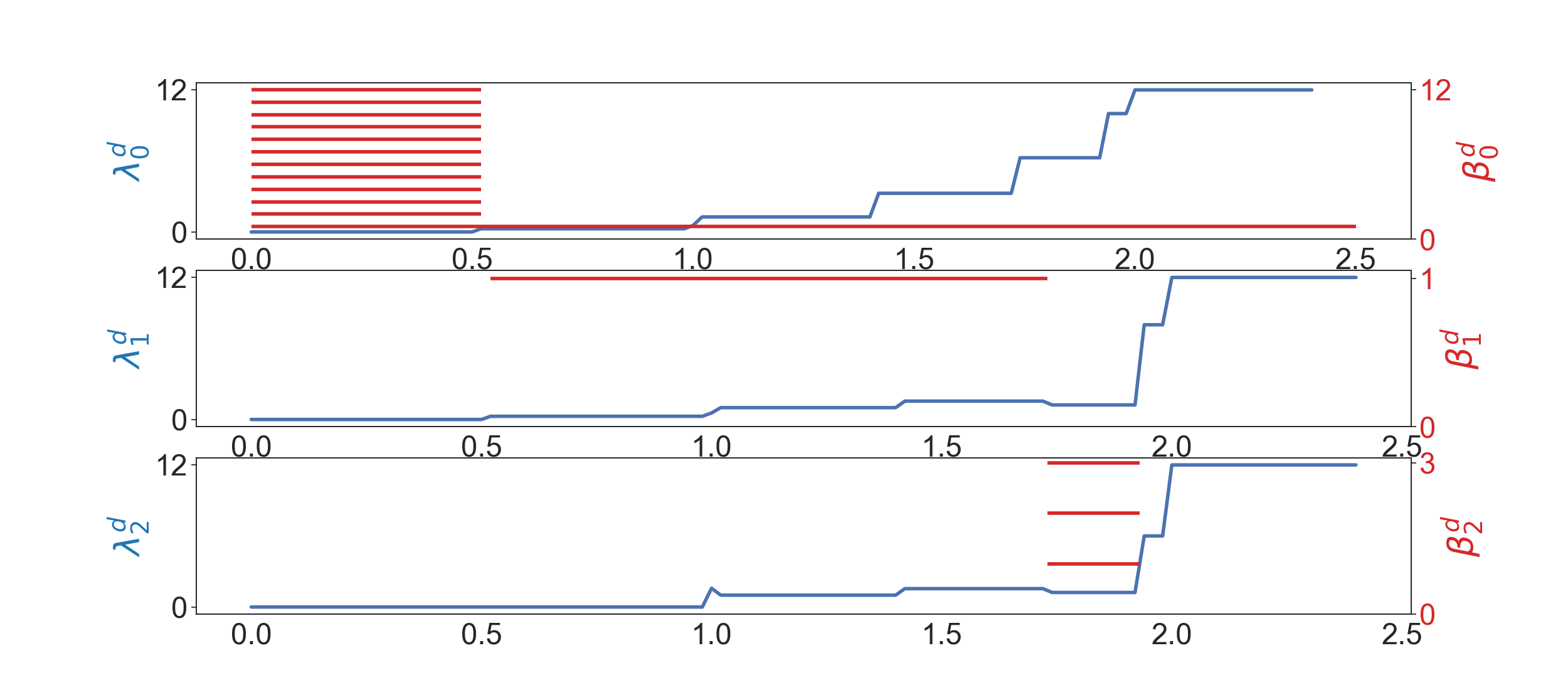}
            %\caption{{\color{red}The results for the Vietoris-Rips filtration constructed from the 8 points.}}
            %\caption{}
            \label{}
        \end{subfigure}
        \caption{}
        \label{12gon}
    \end{figure}
    We consider the Vietoris-Rips filtration of the 12 vertices of the regular 12-gon (Figure \ref{12gon}).
    %The 0-dimensional Betti numbers only tells us that when $d\approx 0.5$, all nodes are connected. 
    %The minimal nonzero eigenvalues changes at different $d$, indicating the forming of new edges.
    We see that $\lambda^d_2$ goes up before any homology class is born, which means that $2$-simplices form but no $2$-dimensional holes are born yet.
    This phenomenon can be observed in many results of \cite{wang2021hermes}.
\end{example}

\subsection{Matrix representations of a (persistent) Laplacian}

%The calculation of matrix representations of a combinatorial Laplacian is simple. 
Since $C_q$ has a canonical orthonormal basis,
the matrix representation of $\partial_q^{\ast}$ is the transpose of the matrix representation of $\partial_q$.
In the computation of persistent Laplacians, the difficult part is the calculation of the up persistent Laplacian, because we need to determine $C_{q+1}^{X,Y}$, which may not have a canonical orthonormal basis.
We can obtain a basis of $C_{q+1}^{X,Y}$ by performing a column reduction for $\partial_{q+1}^Y$, 
or directly calculate the matrix representation of the up persistent Laplacian by Schur complement \cite{memoli2022persistent}.
We give two examples below.

\begin{example}
    \begin{figure}[htbp]
        \centering
        \begin{subfigure}{0.30\textwidth}
          \centering
          \begin{tikzpicture}[decoration={markings, mark= at position 0.75 with {\arrow{triangle 45}[scale=4]}}]      
              \coordinate (0) at (0,0);
              \coordinate (1) at (-1.732, 1);
              \coordinate (2) at (-1.732, -1);
              \coordinate (3) at (0, -2);
              \node at (0.25,0) {$0$};
              \node at (-2,1) {$1$};
              \node at (-2,-1) {$2$};
              \node at (0.25, -2) {$3$};
              \draw (0) -- (1) -- (2) -- (0);
              \draw[fill=black] (0) circle (2pt);
              \draw[fill=black] (1) circle (2pt);
              \draw[fill=black] (2) circle (2pt);
              \draw[fill=black] (3) circle (2pt);
          \end{tikzpicture}
          \caption{}
        \end{subfigure}
        \begin{subfigure}{.30\textwidth}
          \centering
          \begin{tikzpicture}[decoration={markings, mark= at position 0.75 with {\arrow{triangle 45}[scale=4]}}]      
            \coordinate (0) at (0,0);
            \coordinate (1) at (-1.732, 1);
            \coordinate (2) at (-1.732, -1);
            \coordinate (3) at (0, -2);
            \fill[black!20!white] (0) -- (1) -- (2) -- cycle;
            \fill[black!20!white] (0) -- (2) -- (3) -- cycle;
            \draw (0) -- (1) -- (2) -- (0) -- (3) -- (2);
            \node at (0.25,0) {$0$};
            \node at (-2,1) {$1$};
            \node at (-2,-1) {$2$};
            \node at (0.25, -2) {$3$};
            \draw[fill=black] (0) circle (2pt);
            \draw[fill=black] (1) circle (2pt);
            \draw[fill=black] (2) circle (2pt);
            \draw[fill=black] (3) circle (2pt);
        \end{tikzpicture}
        \caption{}
        \end{subfigure}
      \caption{$X = \{0,1,2,3,01,12,02\}$ and $Y = \{0,1,2,3,01,12,23,03,02,012,023\}$.}
      \label{PLsimple case}
    \end{figure}   
    When $C_{q+1}^{X,Y}$ is generated by some $(q+1)$-simplices in $Y$, the calculation of $\partial_q^{X,Y}$ is relatively easy.
    For $X$ and $Y$ shown in Figure \ref{PLsimple case}, we compute $\partial_2^{X,Y}$.
    The matrix representation of $\partial_2^Y$ is 
    \begin{align*}
        \bordermatrix{
            ~  & 012  & 023 \cr
            01 & 1 & 0  \cr
            12 & 1 & 0 \cr
            02 & -1 & 1  \cr
            23 & 0 & 1 \cr
            03 & 0 & -1 \cr
        },
    \end{align*}
    so $C_2^{X,Y}$ is generated by $012$,
    then the matrix representation of $\partial_2^{X, Y}$ is 
    \begin{align*}
        \bordermatrix{
            ~  & 012   \cr
            01 & 1  \cr
            12 & 1  \cr
            02 & -1   \cr
        }.
    \end{align*}  
\end{example}

\begin{example}
    We compute $\Delta_1^{X,Y}$ for $X$ and $Y$ shown in Figure \ref{case2}.
    \begin{figure}[htbp]
        \centering
        \begin{subfigure}{0.30\textwidth}
          \centering
          \begin{tikzpicture}[decoration={markings, mark= at position 0.75 with {\arrow{triangle 45}[scale=4]}}]      
              \coordinate (0) at (0,0);
              \coordinate (1) at (-1.732, 1);
              \coordinate (2) at (-1.732, -1);
              \coordinate (3) at (0, -2);
              \draw  (0) -- (1);
              \draw  (1) -- (2); 
              \draw  (0) -- (3); 
              \draw  (2) -- (3); 
              \node at (0.25,0) {$0$};
              \node at (-2,1) {$1$};
              \node at (-2,-1) {$2$};
              \node at (0.25, -2) {$3$};
              \draw[fill=black] (0) circle (2pt);
              \draw[fill=black] (1) circle (2pt);
              \draw[fill=black] (2) circle (2pt);
              \draw[fill=black] (3) circle (2pt);
          \end{tikzpicture}
          \caption{}
        \end{subfigure}
        \begin{subfigure}{.30\textwidth}
          \centering
          \begin{tikzpicture}[decoration={markings, mark= at position 0.75 with {\arrow{triangle 45}[scale=4]}}]      
            \coordinate (0) at (0,0);
            \coordinate (1) at (-1.732, 1);
            \coordinate (2) at (-1.732, -1);
            \coordinate (3) at (0, -2);
            \fill[black!20!white] (0) -- (1) -- (2) -- cycle;
            \fill[black!20!white] (0) -- (2) -- (3) -- cycle;
            \draw  (0) -- (1);
            \draw  (1) -- (2);
            \draw  (0) -- (2); 
            \draw  (2) -- (3);
            \draw (0) -- (3);
            \node at (0.25,0) {$0$};
            \node at (-2,1) {$1$};
            \node at (-2,-1) {$2$};
            \node at (0.25, -2) {$3$};
            \draw[fill=black] (0) circle (2pt);
            \draw[fill=black] (1) circle (2pt);
            \draw[fill=black] (2) circle (2pt);
            \draw[fill=black] (3) circle (2pt);
        \end{tikzpicture}
          \caption{}
        \end{subfigure}
        \caption{$X= \{0,1,2,3,01,12,23,03\}$ and $Y = \{0,1,2,3,01,12,23,03,02,012,023\}$.}
        \label{case2}
    \end{figure}
    The matrix representation of $\partial_2^Y$ is
    \begin{align*}
        \bordermatrix{
            ~  & 012  & 023 \cr
            01 & 1 & 0  \cr
            12 & 1 & 0 \cr
            23 & 0 & 1 \cr
            03 & 0 & -1 \cr  
            02 & -1 & 1  \cr
        }.
    \end{align*}
    Our goal is to make the submatrix 
    \begin{align*}
        \bordermatrix{
            ~  & 012  & 023 \cr
            02 & -1 & 1  \cr
        }
    \end{align*}in column echelon form.
    We apply one column reduction and get
    \begin{align*}
        \bordermatrix{
            ~  & 012  & 023+012 \cr
            01 & 1 & 1  \cr
            12 & 1 & 1 \cr
            23 & 0 & 1 \cr
            03 & 0 & -1 \cr 
            02 & -1 & 0  \cr
        }.   
    \end{align*}
    Therefore, $C_2^{X,Y}=\text{span}(023+012)$ and one matrix representation of $\partial_2^{X,Y}$ is 
    \begin{align*}
        \bordermatrix{
            ~  & 023+012 \cr
            01 &  1  \cr
            12 &  1 \cr
            23 &  1 \cr
            03 &  -1 \cr 
        }.
    \end{align*}
    For any two spaces $V$ and $W$ and $f:V \to W$, if we choose arbitrary bases of $V$ and $W$ and 
    take a matrix representation $M_f$ of $f$, then the matrix representation $M_{f^{\ast}}$ of $f^{\ast}$ is
    $P^{-1}M_f^T Q$ where $P$ and $Q$ are inner product matrices of $V$ and $W$, respectively.
    If we use $\{023+012\}$ as the basis of $C_2^{X,Y}$, then the inner product matrix of $C_2^{X,Y}$ is $2$ (the norm of $023+012$).
    The corresponding matrix representation of $(\partial_2^{X,Y})^{\ast}$ is 
    \begin{align*}
        \frac{1}{2}\begin{pmatrix}
            1 & 1 & 1 & -1 \\
        \end{pmatrix}
    \end{align*} and the matrix representation of the up persistent Laplacian is 
    \begin{align*}
        \frac{1}{2}\bordermatrix{
            ~  & 01  & 12 & 23 & 03\cr
            01 & 1 & 1 & 1 & -1  \cr
            12 & 1 & 1 & 1 & -1 \cr
            23 & 1 & 1 & 1 & -1\cr
            03 & -1 & -1 & -1 & 1\cr 
        }.
    \end{align*}
    Another way to compute the up persistent Laplacian is as follows.
    We first compute the up Laplacian $\Delta_1(Y)$
    \begin{align*}
        \bordermatrix{
            ~  & 01  & 12 & 23 & 03 & 02\cr
            01 & 1 & 1 & 0 & 0  & -1\cr
            12 & 1 & 1 & 0 & 0 & -1\cr
            23 & 0 & 0 & 1 & -1 & 1\cr
            03 & 0 & 0 & -1 & 1 & -1\cr 
            02 & -1 & -1 & 1 & -1 & 2\cr
        },
    \end{align*}then treat this matrix as a block matrix 
    \begin{align*}
        \bordermatrix{
            ~  & C_1(X)  & 02 \cr
            C_1(X) & A & B \cr
            02 & C & D \cr
        },
    \end{align*}
    and compute the Schur complement $A-BD^{-1}C$.
\end{example}

\subsection{Eigenvectors of a Laplacian}

There are some results concerning the relationship between the spectra of Laplacians and the shape of a simplicial complex \cite{goldberg2002combinatorial,horak2013spectra}.
How do we interpret eigenvectors of a Laplacian?
For an eigenvector of a $q$-th combinatorial Laplacian, 
we can look at the shape of $q$-simplices where the eigenvector has support 
(signs are arbitrary because they are affected by the fixed ordering of vertices).
Empirical observations \cite{krishnagopal2021spectral, muhammad2006control, wei2022hodge} suggest that:
(a) harmonic eigenvectors (eigenvectors of zero eigenvalues) have support near $q$-dimensional ``holes'' (or vertices in a connected component when $q=0$);
%Ebli and Spreemann \cite{ebli2019notion} applied spectral clustering to harmonic eigenvectors. 
%For example, for a simplicial complex that has two ``holes'', this method can separate the simplicial complex into two parts where each part has one ``hole''.
(b) nonharmonic eigenvectors (eigenvectors of nonzero eigenvalues) have support near ``clusters'' of $q$-simplices. 
As to persistent Laplacians, very little is known about the topological interpretation of eigenvalues and eigenvectors.

\section{Generalizations of (persistent) Laplacians}
 
From a theoretical point of view, it is natural to ask whether a persistent Laplacian can be defined in other settings such that persistent Hodge theorem still holds. 
From a practical point of view, generalizations of persistent Laplacians motivated by the need to integrate non-geometrical data that are important for specific problems.
In the next few sections we will introduce structures such as cellular (co)sheaves, digraphs, and hyper(di)graphs, 
and then review some recent advances in their homology and Laplacians. 
We will also discuss Dirac operators and $N$-chain complexes at the end of this section.
%{\color{} think carefully how to talk about this}
%To address this issue, bio-inspired element-specific persistent homology was proposed. 
%The idea is to split the input point cloud into multiple subgroups to better characterize physical and chemical interactions \cite{cang2017topologynet, cang2018integration}. 
%Another more elegant approach is persistent cohomology, which allows the incorporation of non-geometrical information into the topological invariants \cite{cang2020persistent}. 
%Sheaf Laplacians and persistent sheaf Laplacians \cite{wei2021persistent} offer an alternative approach to integrate non-geometric information. 
 
\subsection{Differential graded inner product spaces}

It has been noted early that
persistent Laplacians can be defined analogously for differential graded inner product spaces
and persistent Hodge theorem can be proved.
A \emph{differential graded inner product space} $(V,d)$ is just a chain complex 
\[
    \begin{tikzcd}[column sep = large]
        \centering
    \cdots \arrow[r, "d_{q+2}"] & 
    V_{q+1} \arrow[r, "d_{q+1}"] & 
    V_{q}  \arrow[r, "d_q"] & 
    V_{q-1} \arrow[r, "d_{q-1}"] & 
    \cdots
    \end{tikzcd} 
\]
whose chain groups are inner product spaces.  When we say $(V,d^{V})$ is a subspace of $(W,d^{W})$, 
we mean that the inner space structure and boundary operator $d^V$ of $(V,d^{V})$ are inherited from $(W,d^{W})$. 
For a pair of differential graded inner product spaces $(V,d^{V}) \subset (W, d^{W})$, the $q$-th persistent homology group is defined analogously by 
\begin{align*}
    \iota^{\bullet} (H_q(V)) \cong \frac{\ker d_q^{V}}{\ker d_q^{V} \cap \im d_{q+1}^{W}}.
\end{align*}
Observe that $\ker d_q^{V} \cap \im d_{q+1}^{W} = V_q \cap \im d_{q+1}^{W}$. The preimage of $V_q \cap \im d_{q+1}^{W}$ under $d_{q+1}^W$ is just $(d_{q+1}^W)^{-1}(V_q) = \{w \in W_{q+1} \mid d_{q+1}^Ww \in V_q\}$. 
Hence, $\ker d_q^V \cap \im d_{q+1}^W$ is the image of $\pi d_{q+1}^W \vert_{(d_{q+1}^W)^{-1}(V_q)}: (d_{q+1}^W)^{-1}(V_q) \to V_q$, where $\pi = \iota^{\dag}$, the projection map from $W$ to $V$. 
We denote $\pi d_{q+1}^W\vert_{(d_{q+1}^W)^{-1}(V_q)}$ by $d_{q+1}^{V,W}$, and $(d_{q+1}^W)^{-1}(V_q)$ by $\Theta_{q+1}^{V,W}$. 
These maps are shown in the following diagram
\[
\begin{tikzcd}[column sep = large]
        V_{q+1} \arrow[rr, "d_{q+1}^V"] \arrow[dd, hook, dashed] 
         &
          & V_{q} \arrow[ld, "(d_{q+1}^{V,W})^{\ast}", shift left] \arrow[rr, "d_q^V", shift left] \arrow[dd, hook, dashed] 
           &
            & V_{q-1} \arrow[ll, "(d_q)^{\ast}", shift left] 
            \\
         & \Theta_{q+1}^{V,W} \arrow[ld, hook, dashed] \arrow[ru, "d_{q+1}^{V,W}", shift left] 
          & 
           &
            &\\
        W_{q+1} \arrow[rr, "d_{q+1}^W"] 
         & 
          & W_q 
           &
            & 
    \end{tikzcd}
\]
where hooked dashed arrows represent inclusion maps.
We define the \emph{$q$-th persistent Laplacian} $\Delta_{q}^{V,W}: V_q \to V_q$ by 
\begin{align*}
    (d_q^V)^{\ast}d_q^V + d_{q+1}^{V,W}(d_{q+1}^{V,W})^{\ast}.
\end{align*}
Since $d_q^Vd_{q+1}^{V,W} =0$, we can prove persistent Hodge theorem 
\begin{align*}
    \ker \Delta_{q}^{V,W} \cong \frac{\ker d_q^V}{\ker d_q^V \cap \im d_{q+1}^W}
\end{align*}
in a similar manner. 
Many generalizations of persistent Laplacians implicitly use this formulation.
Liu et al. \cite{liu2023algebraic} first defined persistent Laplacians in the setting of differential graded inner product spaces and showed how to construct a persistent Laplacian for an inner product preserving chain map.

\subsection{Persistent Laplacians for simplicial maps}

The classical filtration of simplicial complexes only represents one type of shape evolution.
We also need tools to study more general shape evolution, such as the sparsification of a simplicial complex.
This requires us to consider general simplicial maps rather than inclusion maps.
G{\"u}len et al. \cite{gulen2023generalization} developed a theory of persistent Laplacians for a simplicial map.
Suppose $f: X \to Y$ is a simplicial map, 
\[
    \begin{tikzcd}
        \cdots \arrow[r]
        & C_{q+1}(X) \arrow[r, "\partial_{q+1}^X"] \arrow[d, "f_{q+1}"] 
          & C_{q}(X) \arrow[r, "\partial_q^X"] \arrow[d, "f_q"] 
            & C_{q-1}(X) \arrow[r] \arrow[d, "f_{q-1}"]
              & \cdots
            \\
        \cdots \arrow[r]
        & C_{q+1}(Y) \arrow[r, "\partial_{q+1}^Y"] 
          & C_q(Y) \arrow[r, "\partial_q^Y"] 
            & C_{q-1}(Y) \arrow[r]
              & \cdots 
    \end{tikzcd}    
\]where $f_q: C_q(X)\to C_q(Y)$ is induced by $f$.
Different from the original $q$-th persistent Laplacian for an inclusion map,
we need to define two subspaces 
\begin{align*}
    C_{q+1}(Y)\supset C_{q+1}^{Y \leftarrow X} = \{c\in C_{q+1}(Y)\mid \partial_{q+1}^Y(c) \in f_q(\ker \partial_q^X)\}
\end{align*}
and 
\begin{align*}
    C_{q-1}(X)\supset C_{q-1}^{X\to Y} = \{c\in C_{q-1}^X\mid (\partial_q^X)^{\ast}(c) \in (\ker f_q)^{\bot}\}
\end{align*}
and then use the restrictions of $\partial_{q+1}^Y$ and $(\partial_q^X)^{\ast}$ to them to construct the $q$-th persistent Laplacian for $f$.
The $q$-th persistent Laplacian for a simplicial map has a more symmetric expression,
and the proof of persistent Hodge theorem is more complicated. 

\subsection{Weighted simplicial complexes}

A simplicial complex whose simplices have weights is generally called a \emph{weighted simplicial complex}.
The weights can be geometrical, such as angles between simplices, volumes of simplices, or non-geometrical such as numbers of scientific papers coauthored by groups of people.
Many theories and models involving weighted simplicial complexes exist (e.g., \cite{baccini2022weighted, battiloro2023topological, chung1996combinatorial, courtney2017weighted, petri2013topological, sharma2017weighted}).
Here we focus on the theory of weighted simplicial complexes proposed by Robert J. MacG. Dawson \cite{dawson1990homology} and 
later developed in \cite{bura2021weighted, bura2022computational, li2022weighted, ren2021weighted, ren2018weighted, wu2018weighted,wu2020discrete}.
A weighted simplicial complex is a simplicial complex where each simplex $\sigma$ has a weight $w(\sigma)$ valued in a commutative ring $R$, 
such that if $\sigma \leqslant \tau$, then $w(\tau)$ is divisible by $w(\sigma)$.
The \emph{weighted chain complex} of a weighted simplicial complex $X$ is defined as follows. 
Let $C_q(X, w)$ be the set of formal sums of $q$-simplices with coefficients in $R$ (if $w(\sigma)$ is zero then we do not include $\sigma$ in any formal sum).
For $\sigma = [v_{a_0},\dots,v_{a_q}]$, we denote the face $[v_{a_0},\dots,\hat{v}_{a_i},\dots, v_{a_q}]$ by $d_i \sigma$.
The \emph{weighted boundary operator} $\partial$ is given by
\begin{align*}
    \partial(\sigma) = \sum_{i=0}^q \frac{w(\sigma)}{w(d_i \sigma)} (-1)^i d_i\sigma.
\end{align*}
As $w(\sigma)$ is divisible by $w(d_i \sigma)$, the weighted boundary operator is well-defined.
%When $w$ is constant, the weighted boundary map is the same as the boundary map of a simplicial complex.
We still have $\partial^2=0$, because for $0\leq i < j \leq q$,
\begin{align*}
    \frac{w(\sigma)}{w(d_i\sigma)} \frac{w(d_i\sigma)}{w(d_{j-1}d_i\sigma)} = \frac{w(\sigma)}{w(d_j \sigma)} \frac{w(d_j\sigma)}{w(d_id_j\sigma)} = \frac{w(\sigma)}{w(d_id_j\sigma)}.
\end{align*}
Therefore, weighted homology groups can be defined analogously.
Wu et al. \cite{wu2018weighted} pointed out that in the proof of $\partial^2=0$, what really matters is the quotient of weights. 
If we write $w(\tau)/w(\sigma)$ as $\phi(\tau, \sigma)$, 
then the equality
\begin{align*}
    \frac{w(\sigma)}{w(d_i\sigma)} \frac{w(d_i\sigma)}{w(d_{j-1}d_i\sigma)} = \frac{w(\sigma)}{w(d_j \sigma)} \frac{w(d_j\sigma)}{w(d_id_j\sigma)} 
\end{align*}
becomes 
\begin{align*}
    \phi(d_i\sigma, d_{j-1}d_i \sigma)\phi(\sigma,d_i\sigma) = \phi(d_j\sigma, d_id_j \sigma)\phi(\sigma, d_j\sigma),
\end{align*}
which means that any $\phi: X \times X \to R$ satisfying this equality induces a 
($\phi$-weighted) boundary operator
\begin{align*}
    \partial_q(\sigma) = \sum_{i=0}^q (-1)^i\phi(\sigma, d_i \sigma)d_i \sigma.
\end{align*}
A simplicial complex paired with a generalized weight function $\phi$ is called a $\phi$-weighted simplicial complex. 

\begin{example}
    \cite{wu2018weighted} A weighted polygon is a polygon with $\phi(\{v_i,v_j\}, v_i) = \alpha_i \in \intg$ (Figure \ref{weightedPolygon}).
    \begin{figure}[htbp]
        \centering
          \begin{tikzpicture}[decoration={markings, mark= at position 0.75 with {\arrow{triangle 45}[scale=4]}}]      
              \coordinate (0) at (1,0);
              \coordinate (1) at (2, 2);
              \coordinate (2) at (0, 4);
              \coordinate (3) at (-2, 2);
              \coordinate (4) at (-1,0);
              \draw  (0) -- (1) -- (2) -- (3) -- (4) -- (0);
              \draw pic ["$\alpha_0$", draw,thick, angle radius=0.4cm, angle eccentricity=1.7] {angle = 1--0--4};
              \draw pic ["$\alpha_1$", draw,thick, angle radius=0.4cm, angle eccentricity=1.7] {angle = 2--1--0};
              \draw pic ["$\alpha_2$", draw,thick, angle radius=0.4cm, angle eccentricity=1.7] {angle = 3--2--1};
              \draw pic ["$\alpha_3$", draw,thick, angle radius=0.4cm, angle eccentricity=1.7] {angle = 4--3--2};
              \draw pic ["$\alpha_4$", draw,thick, angle radius=0.4cm, angle eccentricity=1.7] {angle = 0--4--3};
              \node (A) [right= of 0]  {} ;
              \draw[fill=black] (0) circle (2pt);
              \draw[fill=black] (1) circle (2pt);
              \draw[fill=black] (2) circle (2pt);
              \draw[fill=black] (3) circle (2pt);
              \draw[fill=black] (4) circle (2pt);
          \end{tikzpicture}
          \caption{A weighted polygon.}
          \label{weightedPolygon}
    \end{figure}
    The matrix representation of $\partial_1$ is
    \begin{align*}
        \bordermatrix{
            ~  & v_0v_1  & v_0v_4 & v_1v_2 & v_2v_3 & v_3v_4\cr
            v_0 & -\alpha_0 & -\alpha_0 & 0 & 0 & 0  \cr
            v_1 & \alpha_1 & 0 & -\alpha_1 & 0 & 0 \cr
            v_2 & 0 & 0 & \alpha_2 & -\alpha_2 & 0 \cr
            v_3 & 0 & 0 & 0 & \alpha_3 & -\alpha_3 \cr 
            v_4 & 0 & \alpha_4 & 0 & 0 & \alpha_4 \cr
        }
    \end{align*}
    and the resulting weighted $H_0$ is dependent on $\alpha_i$. 
    The weighted homology of weighted polygons might be useful for analyzing ring structures in biomolecules.
\end{example} 

We have emphasized that a point cloud can be studied by building a filtration of simplicial complexes.
If we want to distinguish some points from other points, we can assign weights and building a filtration of weighted simplicial complexes \cite{ren2018weighted}.
We may also consider weighted versions of (persistent) Laplacians \cite{wu2018weighted}.

\begin{example}
    Suppose each point $v$ in a point cloud has weight $w(v)$.
    We can associate any simplex $\{v_{a_0},\dots, v_{a_q}\}$ the product weight \cite{ren2018weighted}    
    \begin{align*}
        \prod_{i=0}^q w(v_{a_i}).
    \end{align*}
    Since the weighted boundary map can be given by 
    \begin{align*}
        \partial(\sigma) = \sum_{i=0}^q w(v_{a_i}) (-1)^i d_i\sigma.
    \end{align*}
    We can just define the $q$-th chain group as the space generated by $q$-simplices without worrying about zero weights.
\end{example}

\begin{example} 
    Suppose a point cloud contains two types of points $A$ and $B$.
    We can assign weights $\{0,1\}$ to $\{A,B\}$, and compute weighted homology and Laplacians using product weighting.   
    At least when a point cloud is simple, 
    weighted combinatorial Laplacians can be used to differentiate among different patterns of distribution of $A$ and $B$.
    For a point cloud of four points $\{(0,0), (1,0), (1,1), (0,1)\}$ there are five configurations (shown in Figure \ref{conf}) that include at least one point whose weight is 1.
    Results of weighted Laplacians are shown in Figure \ref{confRes}.
    \begin{figure}[htbp]
        \centering
        \begin{subfigure}{0.18\textwidth}
          \centering
          \begin{tikzpicture}[decoration={markings, mark= at position 0.75 with {\arrow{triangle 45}[scale=4]}}]      
              \coordinate (0) at (0,0);
              \coordinate (1) at (1, 0);
              \coordinate (2) at (1, 1);
              \coordinate (3) at (0, 1);
              %\draw  (0) -- (1) -- (2) -- (3) -- (0);
              \draw[fill=black] (0) circle (2pt);
              \draw[fill=black] (1) circle (2pt);
              \draw[fill=black] (2) circle (2pt);
              \draw[fill=black] (3) circle (2pt);
          \end{tikzpicture}
          \caption{}
        \end{subfigure}
        \begin{subfigure}{0.18\textwidth}
            \centering
            \begin{tikzpicture}[decoration={markings, mark= at position 0.75 with {\arrow{triangle 45}[scale=4]}}]      
                \coordinate (0) at (0,0);
                \coordinate (1) at (1, 0);
                \coordinate (2) at (1, 1);
                \coordinate (3) at (0, 1);
                %\draw  (0) -- (1) -- (2) -- (3) -- (0);
                \draw[fill=black] (0) circle (2pt);
                \draw[fill=black] (1) circle (2pt);
                \draw[fill=black] (2) circle (2pt);
                \draw[fill=white] (3) circle (2pt);
            \end{tikzpicture}
            \caption{}
          \end{subfigure}
          \begin{subfigure}{0.18\textwidth}
            \centering
            \begin{tikzpicture}[decoration={markings, mark= at position 0.75 with {\arrow{triangle 45}[scale=4]}}]      
                \coordinate (0) at (0,0);
                \coordinate (1) at (1, 0);
                \coordinate (2) at (1, 1);
                \coordinate (3) at (0, 1);
                %\draw  (0) -- (1) -- (2) -- (3) -- (0);
                \draw[fill=black] (0) circle (2pt);
                \draw[fill=black] (1) circle (2pt);
                \draw[fill=white] (2) circle (2pt);
                \draw[fill=white] (3) circle (2pt);
            \end{tikzpicture}
            \caption{}
          \end{subfigure}
          \begin{subfigure}{0.18\textwidth}
            \centering
            \begin{tikzpicture}[decoration={markings, mark= at position 0.75 with {\arrow{triangle 45}[scale=4]}}]      
                \coordinate (0) at (0,0);
                \coordinate (1) at (1, 0);
                \coordinate (2) at (1, 1);
                \coordinate (3) at (0, 1);
                %\draw  (0) -- (1) -- (2) -- (3) -- (0);
                \draw[fill=black] (0) circle (2pt);
                \draw[fill=white] (1) circle (2pt);
                \draw[fill=black] (2) circle (2pt);
                \draw[fill=white] (3) circle (2pt);
            \end{tikzpicture}
            \caption{}
          \end{subfigure}
          \begin{subfigure}{0.18\textwidth}
            \centering
            \begin{tikzpicture}[decoration={markings, mark= at position 0.75 with {\arrow{triangle 45}[scale=4]}}]      
                \coordinate (0) at (0,0);
                \coordinate (1) at (1, 0);
                \coordinate (2) at (1, 1);
                \coordinate (3) at (0, 1);
                %\draw  (0) -- (1) -- (2) -- (3) -- (0);
                \draw[fill=black] (0) circle (2pt);
                \draw[fill=white] (1) circle (2pt);
                \draw[fill=white] (2) circle (2pt);
                \draw[fill=white] (3) circle (2pt);
            \end{tikzpicture}
            \caption{}
          \end{subfigure}
          \caption{}
          \label{conf}
    \end{figure}    
    \begin{figure}[htbp]
        \centering
        \begin{subfigure}{0.30\textwidth}
            \centering
            \includegraphics[width=0.80\textwidth]{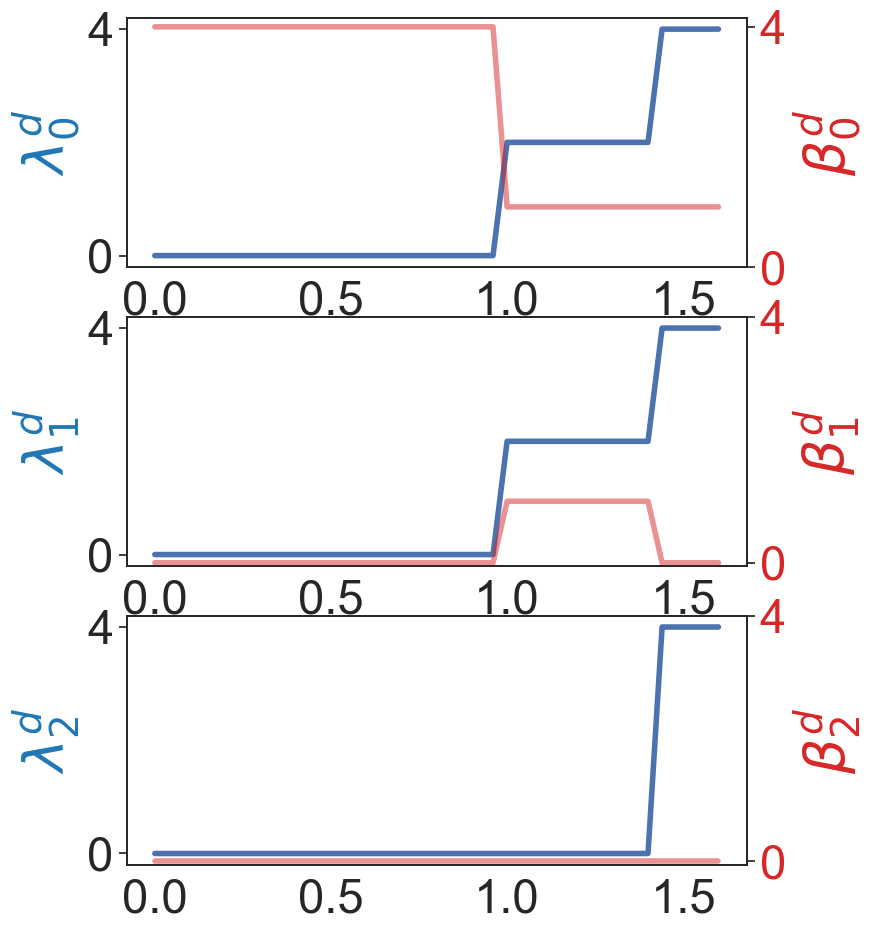}
            \caption{$BBBB$}
        \end{subfigure}
        \begin{subfigure}{0.30\textwidth}
            \centering
            \includegraphics[width=0.80\textwidth]{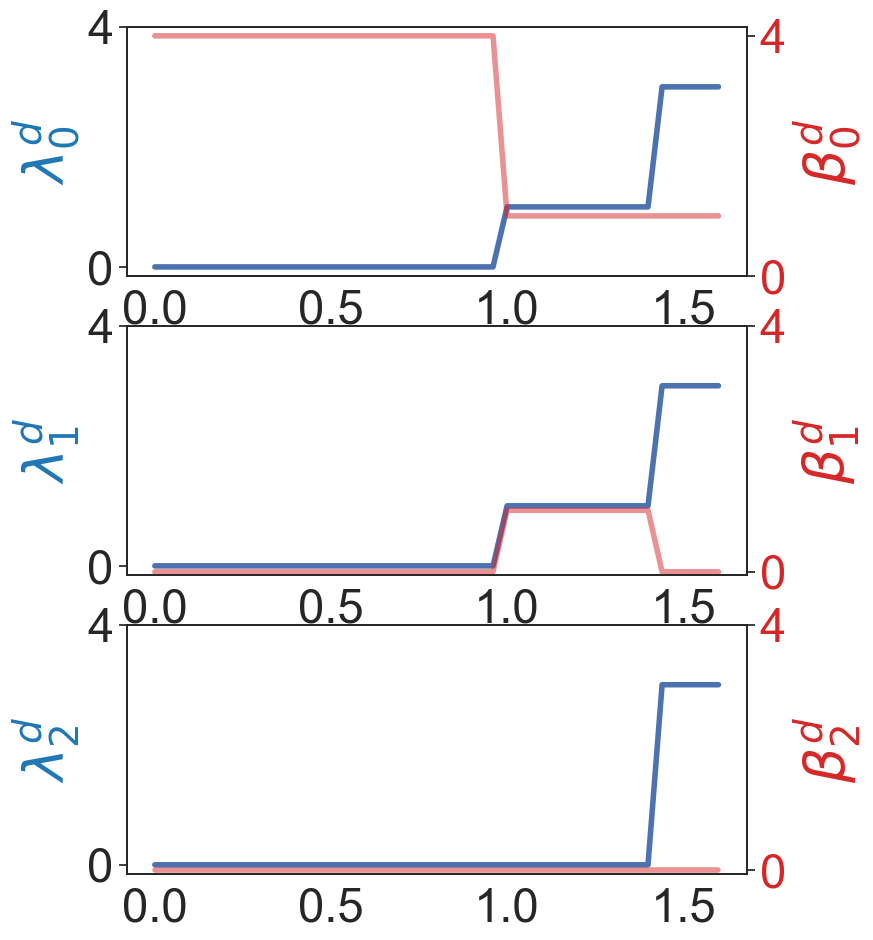}
            \caption{$ABBB$}
          \end{subfigure}
        \begin{subfigure}{0.30\textwidth}
            \centering
            \includegraphics[width=0.80\textwidth]{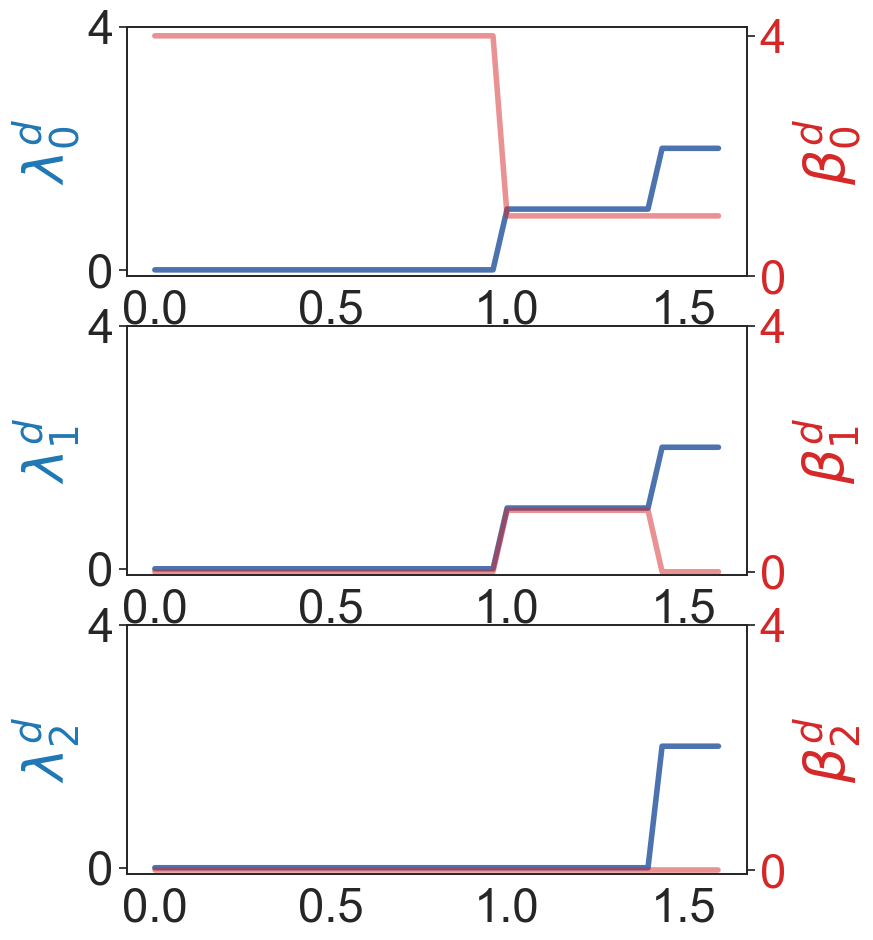}
            \caption{$AABB$}
        \end{subfigure}
        \begin{subfigure}{0.30\textwidth}
            \centering
            \includegraphics[width=0.80\textwidth]{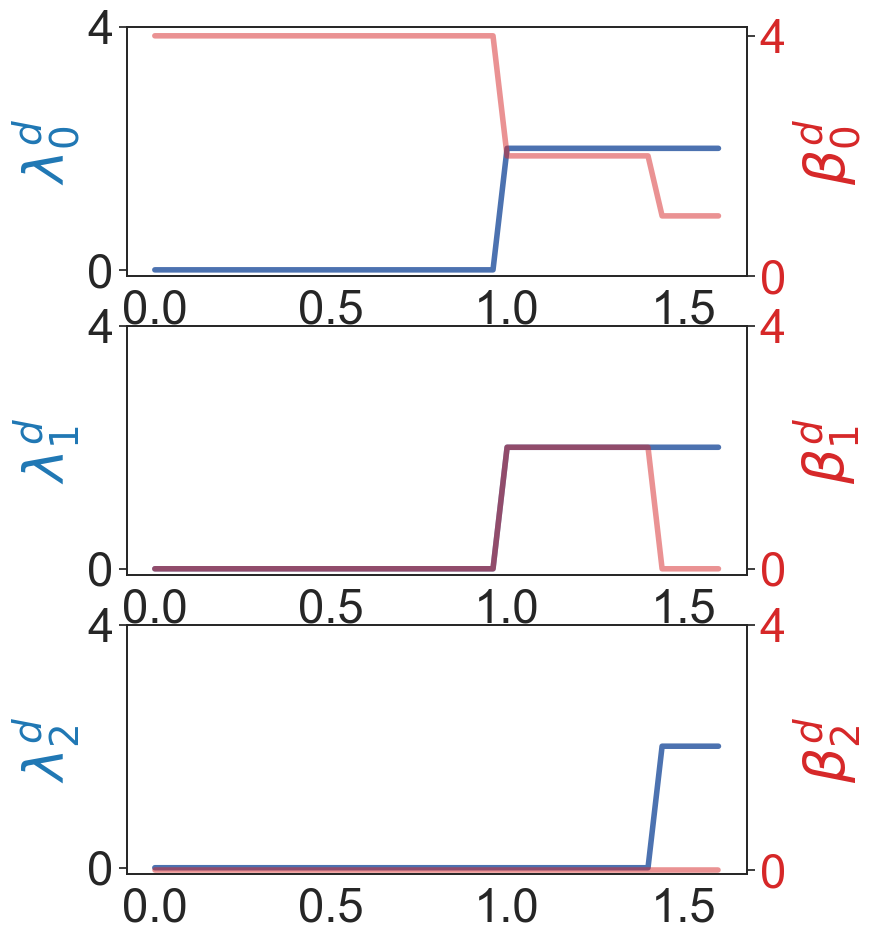}
            \caption{$ABAB$}
        \end{subfigure}
        \begin{subfigure}{0.30\textwidth}
            \centering
            \includegraphics[width=0.80\textwidth]{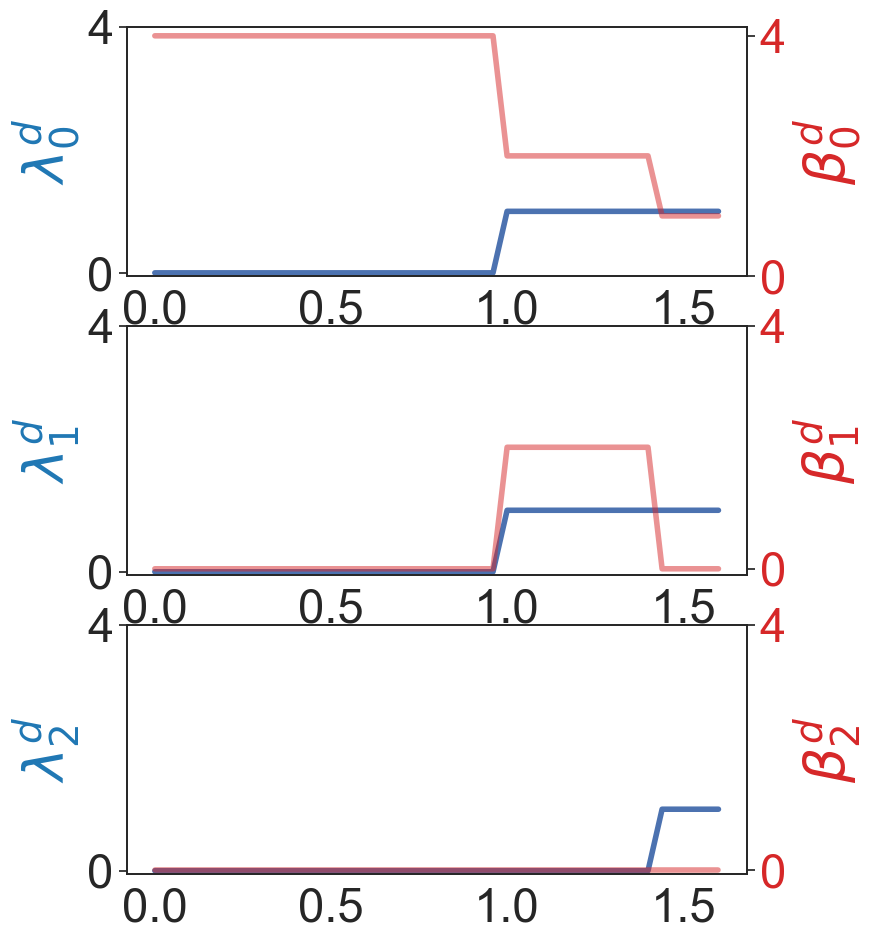}
            \caption{$AAAB$}
        \end{subfigure}
        \caption{$\lambda^d_q$ is the minimal nonzero eigenvalue of the $q$-th weighted combinatorial Laplacian for $X_d$ in a Rips filtration. 
        $\beta^d_q$ is the $q$-th Betti number of $X_d$.}
        \label{confRes}
    \end{figure}
\end{example} 

\subsection{Cellular (co)sheaves}
%To integrate non-spatial information into persistent homology, 
%approaches such as element specific persistent homology \cite{cang2017topologynet,} and persistent cohomology \cite{cang2020persistent} were proposed.
%In this section we introduce an alternative approach.

In a $\phi$-weighted simplicial complex, 
we can imagine that a copy of $R$ resides each simplex and $\phi(\tau, \sigma)$ is a scalar multiplication from the copy on $\tau$ to the copy on $\sigma$ \cite{hansen2019toward}. 
If we associate each simplex with a vector space and designate a linear morphism for every face relation, we will get a \emph{cellular (co)sheaf}.
The theory of cellular (co)sheaves was first introduced in \cite{shepard1985cellular} and 
later gained attention for its application potential (e.g., \cite{curry2014sheaves, hansen2020laplacians, robinson2014topological, robinson2017sheaves, yoon2018cellular}). 
In the recent years, the study of sheaf neural networks has become a trending topic \cite{barbero2022sheaf, bodnar2022neural, gillespie2024bayesian, hansen2020sheaf}. 
Like a weighted simplicial complex, a cellular (co)sheaf is a candidate for modeling complex objects such as molecules.

A cellular cosheaf $\mathscr{F}$ is a simplicial complex $X$ with additional data\footnote{For ease of exposition we have simplified the definition of a cellular (co)sheaf.}.
Each simplex $\sigma$ is assigned a vector space $\mathscr{F}(\sigma)$ (or denoted by $\mathscr{F}_{\sigma}$), referred to as the stalk over $\sigma$, 
and for any face relation $\sigma \leqslant \tau$,  
there is an \emph{extension map} $\mathscr{F}(\sigma \leqslant \tau): \mathscr{F}(\tau) \to \mathscr{F}(\sigma)$ (or denoted by $\mathscr{F}_{\sigma \leqslant \tau}$).
The $q$-th chain group of a cellular cosheaf is the direct sum of stalks over $q$-simplices and
the boundary map $\partial$ is given by 
\begin{align*}
    \partial_q\vert_{\mathscr{F}(\sigma)} = \sum_{i} (-1)^i \mathscr{F}(d_i\sigma \leqslant \sigma).
\end{align*}
The square of this boundary map is 0, if for any face relation $\rho \leqslant \sigma \leqslant \tau$ we have
\begin{align*}
    \mathscr{F}(\rho \leqslant \tau) = \mathscr{F}(\rho \leqslant \sigma) \circ \mathscr{F}(\sigma \leqslant \tau).
\end{align*}
A dual concept is a cellular sheaf. 
For a cellular sheaf $\mathscr{F}$, $\mathscr{F}(\sigma \leqslant \tau)$ is a map from $\mathscr{F}(\sigma)$ to $\mathscr{F}(\tau)$ (called a restriction map).
A sheaf cochain complex can be defined analogously. 
If stalks are inner product spaces, one can equip inner product structures for (co)chain groups.
Applying the construction of combinatorial Laplacians to a cosheaf chain complex or sheaf cochain complex, 
we get (co)sheaf Laplacians \cite{hansen2019toward}. 
It is noted that many sheaves over a digraph only have trivial 0-dimensional cohomology groups \cite{gentleintrosheaves},
but we can still extract some information from sheaf Laplacians.

\begin{example}
    Suppose there is a sheaf $\mathscr{F}$ over the simplicial complex $\{0,1,2,01,02,12\}$,
    then the sheaf coboundary map $\delta^0$ is represented by the block matrix
    \begin{align*}
            \bordermatrix{
                ~               & \mathscr{F}_{0}         & \mathscr{F}_{1}      & \mathscr{F}_{2}\cr
                \mathscr{F}_{01} & -\mathscr{F}_{0\leqslant 01} & \mathscr{F}_{1\leqslant 01} & 0   \cr
                \mathscr{F}_{02} & -\mathscr{F}_{0\leqslant 02} & 0 & \mathscr{F}_{2\leqslant 02}  \cr
                \mathscr{F}_{12} & 0 & -\mathscr{F}_{1\leqslant 12} & \mathscr{F}_{2\leqslant 12}  \cr
            }.
    \end{align*}
    Suppose all stalks are inner product spaces and they are orthogonal to each other, the $0$-th sheaf Laplacian ${\delta^0}^{\ast} \delta^0$ is represented by the block matrix 
    \begin{align*}
        \bordermatrix{
            ~               & \mathscr{F}_0         & \mathscr{F}_1      & \mathscr{F}_2\cr
            \mathscr{F}_0 & \mathscr{F}_{0\leqslant 01}^{\ast}\mathscr{F}_{0\leqslant 01} + \mathscr{F}_{0\leqslant 02}^{\ast}\mathscr{F}_{0\leqslant 02}  & -\mathscr{F}_{0\leqslant 01}^{\ast}\mathscr{F}_{1\leqslant 01}  & -\mathscr{F}_{0\leqslant 02}^{\ast}\mathscr{F}_{2\leqslant 02}   \cr
            \mathscr{F}_1 & -\mathscr{F}_{1\leqslant 01}^{\ast}\mathscr{F}_{0\leqslant 01} & \mathscr{F}_{1\leqslant 01}^{\ast}\mathscr{F}_{1\leqslant 01} + \mathscr{F}_{1\leqslant 12}^{\ast}\mathscr{F}_{1\leqslant 12} & -\mathscr{F}_{1\leqslant 12}^{\ast}\mathscr{F}_{2\leqslant 12}  \cr
            \mathscr{F}_2 & -\mathscr{F}_{2\leqslant 02}^{\ast}\mathscr{F}_{0\leqslant 02} & -\mathscr{F}_{2\leqslant 12}^{\ast}\mathscr{F}_{1\leqslant 12} & \mathscr{F}_{2\leqslant 02}^{\ast}\mathscr{F}_{2\leqslant 02} + \mathscr{F}_{2\leqslant 12}^{\ast}\mathscr{F}_{2\leqslant 12} \cr
        }.        
    \end{align*}
\end{example}
    
Persistent (co)sheaf (co)homology is known to experts \cite{michaelrobinsontutorial, yegnesh2016persistence} 
and a systematical treatment can be found in \cite{russold2022persistent}. 
%Many types of filtration of (co)sheaves can be constructed. 
One type of filtration of sheaves is as follows. 
A sheaf $\mathscr{F}$ over $X$ is a ``subsheaf'' of $\mathscr{G}$ over $Y$
if $X \subset Y$ and stalks and restriction maps of $\mathscr{F}$ are the same as those of $\mathscr{G}$.
To define the $q$-th persistent sheaf Laplacian \cite{wei2021persistent} for $\mathscr{F} \subset \mathscr{G}$, we can first endow chain groups with inner product structures
and dualize everything to make $\mathscr{F}$ and $\mathscr{G}$ cosheaves, such that there is an inclusion chain map between their cosheaf chain complexes.
Then we can define the $q$-th persistent sheaf Laplacian as the $q$-th persistent Laplacian of cosheaf chain complexes.

\subsection{Path homology, flag homology, and digraphs}

The motivation behind path homology is to construct a homology theory of digraphs such that directional information of edges is encoded and higher dimensional homology groups are   non-trivial.
Path homology\footnote{There are other (co)homology theories of digraphs \cite{caputi2023hochschild, lutgehetmann2020computing, masulli2016topology, reimann2017cliques}.} was proposed by Grigor'yan, Lin, Muranov and Yau \cite{grigor2012homologies} and developed in various papers \cite{grigor2014homotopy, grigor2015cohomology, grigor2020path, grigor2017homologies, lin2019weighted}.
A summary of recent advances in path homology of digraphs can be found in \cite{grigor2022advances}.
Recall that a digraph (without self-loops) is a pair $G=(V, E)$ where $E$ is a set of ordered pairs of vertices. 
An \emph{allowed $q$-path} is an ordered finite sequence of vertices $\{x_0, \dots, x_q\}$ such that 
$(x_i, x_{i+1}) \in E$ for all $0\leq i \leq q-1$. 
If we take the space generated by allowed $p$-paths (denoted by $\mathcal{A}_q$) as the $q$-th chain group,
and define the boundary map $\partial_q$ by 
\begin{align*}
    \partial_q \{x_0, \dots, x_q\} = \sum_{i=0}^q (-1)^{i} \{x_0, \dots, \hat{x}_i, \dots, x_q \}
\end{align*}
then formally we can show that $\partial^2=0$.
However, $\partial_q \{x_0, \dots, x_q\}$ may include paths that are not allowed paths. To solve this problem, we need to introduce some general concepts first.

\begin{definition}
    Suppose $X$ is a finite set. An elementary $p$-path is a sequence $[x_0, \dots, x_p]$ of $p+1$ elements of $X$. The space generated by all elementary $p$-paths with coefficient in $\real$ is denoted by $\Lambda_p(X)$. The $q$-th non-regular boundary map is given by 
    \begin{align*}
        \partial_q^{\text{nr}} [x_0, \dots, x_q] = \sum_{i=0}^p [x_0, \dots, \hat{x}_i, \dots, x_q].
    \end{align*}
\end{definition}
One can prove this is a chain complex.
Among all the paths, a path that lingers at a vertex (for some $i$, $x_i = x_{i+1}$) is considered a degenerate path since we are not interested in self-loops. 

\begin{definition}
    A path $[x_0, \dots, x_q]$ over $X$ where $x_i \neq x_{i+1}$ for each $i$ is called regular. The space generated by all regular $q$-paths is denoted by $\mathcal{R}_q$.     
\end{definition}

We define a new boundary operator $\partial_q$ between regular paths. 
When computing $\partial_q([x_0, \dots, x_q])$, we first compute $\partial_q^{\text{nr}}([x_0, \dots, x_q])$ and treat all irregular paths arising from it as zeros. One can still verify that $\partial^2 = 0$ \cite{grigor2019homology}.

Now given a digraph $G=(V, E)$, every $\mathcal{A}_q$ is a subspace of $\mathcal{R}_q$.
\[
    \begin{tikzcd}
    \dots \arrow[r, "\partial_{q+2}"] & 
    \mathcal{R}_{q+1} \arrow[r, "\partial_{q+1}"] & 
    \mathcal{R}_{q}  \arrow[r, "\partial_q"] & 
    \mathcal{R}_{q-1} \arrow[r, "\partial_{q-1}"] & 
    \dots \\ 
     & 
     \mathcal{A}_{q+1} \arrow[u, hook, dashed] &
     \mathcal{A}_q \arrow[u, hook, dashed] &
     \mathcal{A}_{q-1} \arrow[u, hook, dashed] & 
     & %\\  
     %& 
     %D_{q+1}\cap \partial_{q+1}^{-1}(D_q) %\arrow[u, hook, dashed] &
     %D_q\cap \partial_{q}^{-1}(D_{q-1}) \arrow[u, hook, dashed] &
     %D_{q-1}\cap \partial_{q-1}^{-1}(D_{q-2}) \arrow[u, hook, dashed] & 
     %&
    \end{tikzcd}
\]

One way to make $\partial_q: \mathcal{A}_q \to \mathcal{A}_{q-1}$ well-defined is to restrict $\partial_q$ to the subspace $\mathcal{A}_{q} \cap \partial_q^{-1}\mathcal{A}_{q-1}$. 
We have to verify that $\partial_q (\mathcal{A}_{q} \cap \partial_q^{-1}\mathcal{A}_{q-1}) \subset \mathcal{A}_{q-1}\cap \partial_{q-1}^{-1}\mathcal{A}_{q-2}$. $\partial_q (\mathcal{A}_{q} \cap \partial_q^{-1}\mathcal{A}_{q-1}) \subset \mathcal{A}_{q-1}\cap $ is true by definition, and $\partial_q (\mathcal{A}_{q} \cap \partial_q^{-1}\mathcal{A}_{q-1}) \subset \partial_{q-1}^{-1}\mathcal{A}_{q-2}$ is true since $\partial^2 = 0$.
Therefore, we have the chain complex 
\[
    \begin{tikzcd}
    \dots \arrow[r,] & 
    \mathcal{A}_{q+1} \cap \partial_{q+1}^{-1}\mathcal{A}_{q}  \arrow[r, "\partial_{q+1}"] & 
    \mathcal{A}_{q} \cap \partial_q^{-1}\mathcal{A}_{q-1} \arrow[r, "\partial_q"] & 
    \mathcal{A}_{q} \cap \partial_q^{-1}\mathcal{A}_{q-1} \arrow[r] & 
    \dots
    \end{tikzcd}
\]
and the definition of a path homology group is straightforward.
The $q$-th chain group $\mathcal{A}_{q} \cap \partial_q^{-1}\mathcal{A}_{q-1}$ is called the space of $\partial$-invariant $q$-paths on $G$, denoted by $\Omega_q$\footnote{If a digraph is not simple, 
there will be two choices of $\partial_q$ \cite{grigor2012homologies} that might be suitable for different problems \cite{huntsman2020path}.}.
Regarding the geometrical interpretation of path homology,
we only know that non-reduced $H_0$ is the number of connected components of the underlying undirected graph.
%It is not easy to relate higher dimensional path homology groups to features of a digraph. 
Chowdhury et al. \cite{chowdhury2022path} obtained some characterizations of path homologies of certain families of small digraphs. 
Since directional information of edges is encoded in path homology, path homology can be used to distinguish network motifs \cite{chowdhury2018persistent} and isomers in molecular and materials sciences \cite{chen2023path}. 
We can also quantify the significance of a node in a network by observing the changes of path homology when the node is removed \cite{chen2023path}.

Since $\Omega_q$ inherits the inner product structure from $\mathcal{A}_q$, 
the so-called path Laplacians\footnote{Another type of path Laplacians was proposed by Estrada \cite{estrada2012path}
and applied in molecular biology \cite{liu2023persistent}.} can be defined. We can use path Laplacians \cite{gomes2019path, grigor2022advances, wang2023persistent} to distinguish among digraphs that 
have the same path homology.
For example, according to \cite[Theorem 5.4]{grigor2012homologies}, the following two digraphs $G_L$ and $G_R$ (see Figure \ref{GLGR}) have the same path homology. 
But the spectrum of the $0$-th path Laplacian of $G_L$ is $\{0,3,3\}$ and that of $G_R$ is $\{0,2,4,4\}$.

\begin{figure}[htbp]
    \centering
    \begin{subfigure}{0.30\textwidth}
        \centering
        \begin{tikzpicture}[decoration={markings, mark= at position 0.75 with {\arrow{triangle 45}[scale=4]}}, every edge quotes/.style = {auto, font=\footnotesize}]      
            \coordinate (0) at (0,0);
            \coordinate (1) at (0, 2);
            \coordinate (2) at (-1.732, 1);
            \draw [postaction={decorate}] (0) -- (1);
            \draw [postaction={decorate}] (1) -- (2);
            \draw [postaction={decorate}] (2) -- (0); 
            \draw[fill=black] (0) circle (2pt);
            \draw[fill=black] (1) circle (2pt);
            \draw[fill=black] (2) circle (2pt);
        \end{tikzpicture}
        \caption{$G_L$}
        \label{}
    \end{subfigure}
    \begin{subfigure}{0.30\textwidth}
        \centering
        \begin{tikzpicture}[decoration={markings, mark= at position 0.75 with {\arrow{triangle 45}[scale=4]}}]      
            \coordinate (0) at (0,0);
            \coordinate (1) at (0, 2);
            \coordinate (2) at (-1.732, 1);
            \coordinate (3) at (1.732, 1);
            \draw [postaction={decorate}] (0) -- (1);
            \draw [postaction={decorate}] (1) -- (2);
            \draw [postaction={decorate}] (2) -- (0); 
            \draw [postaction={decorate}] (1) -- (3); 
            \draw [postaction={decorate}] (0) -- (3); 
            \draw[fill=black] (0) circle (2pt);
            \draw[fill=black] (1) circle (2pt);
            \draw[fill=black] (2) circle (2pt);
            \draw[fill=black] (3) circle (2pt);
        \end{tikzpicture}
        \caption{$G_R$}
        \label{}
    \end{subfigure}
		\caption{Two digraphs that have the same path homology.}
    \label{GLGR}
\end{figure}
Persistent path homology was proposed by Chowdhury and M{\'e}moli \cite{chowdhury2018persistent} to study a digraph where each edge $e$ has a weight $w(e)$.
A filtration of digraphs $\{G_d\}$ is constructed such that $e \in G_d$ if and only if $w(e) \leq d$.
Wang and Wei \cite{wang2023persistent} introduced persistent path Laplacians and 
demonstrated that persistent path Laplacians can be applied to study molecules, since a lot of information of molecules can be encoded in digraphs.

Flag complexes, also known as clique complexes, is another way to construct homology for digraphs, and they arise naturally in many situations \cite{lutgehetmann2020computing}. 
Jones and Wei \cite{jones2023persistent} introduced persistent directed flag Laplacians as a distinct way of analyzing flag complexes, 
and applied them to analyze protein-ligand binding data. 

\begin{example}
    For a weighted digraph, we can build a filtration $\{G_d\}$ such that $e \in G_d$ iff $w(e) \leq d$.
    Two weighted graphs whose path Betti numbers are the same for every $G_d$ may have different path Laplacians.
    \begin{figure}[htbp]
        \centering
        \begin{subfigure}{0.4\textwidth}
            \centering
            \begin{tikzpicture}
            % Nodes
            \node[circle, draw] (0) at (1,0) {0};
            \node[circle, draw] (1) at (0,1) {1};
            \node[circle, draw] (2) at (0,-1) {2};
            \node[circle, draw] (3) at (-1,0) {3};
            \draw[->,>=stealth] (1) -- node[midway,above] {1} (0);
            \draw[->,>=stealth] (2) -- node[midway,below] {2} (0);
            \draw[->,>=stealth] (1) -- node[midway,right] {0} (2);
            \draw[->,>=stealth] (3) -- node[midway,above] {0} (1);
            \draw[->,>=stealth] (3) -- node[midway,below] {0} (2);
            \end{tikzpicture}
            \caption{}
        \end{subfigure}
        \begin{subfigure}{0.4\textwidth}
            \centering
            \begin{tikzpicture}
                \node[circle, draw] (a) at (1,0) {0};
                \node[circle, draw] (b) at (0,1) {1};
                \node[circle, draw] (c) at (0,-1) {2};
            
                \draw[->,>=stealth] (b) -- node[midway,above] {1} (a);
                \draw[->,>=stealth] (c) -- node[midway,below] {2} (a);
                \draw[->,>=stealth] (b) -- node[midway,right] {0} (c);            
            \end{tikzpicture}   
            \caption{}          
        \end{subfigure}
        \begin{subfigure}{0.4\textwidth}
            \centering
            \includegraphics[width=1\textwidth]{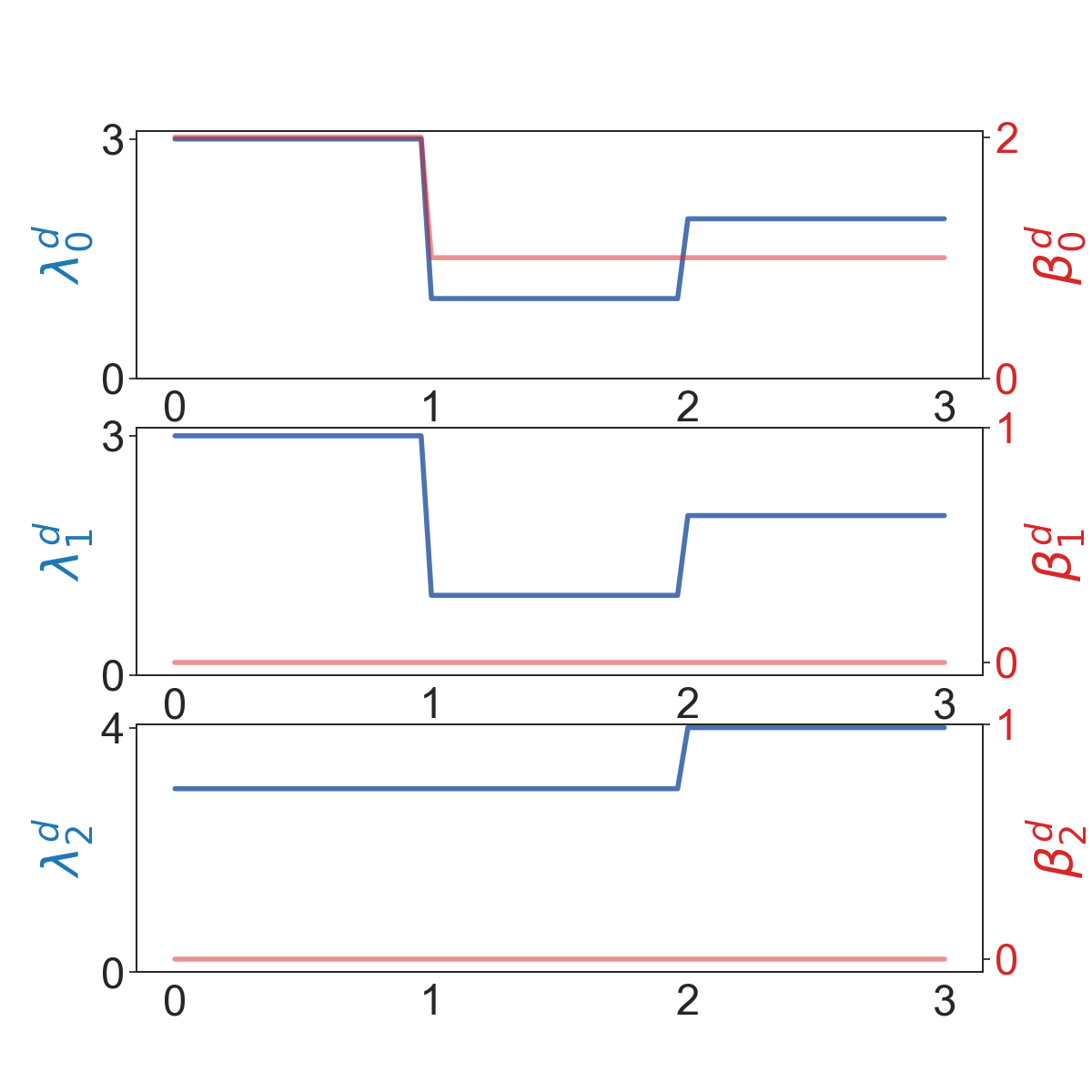}
            \caption{Results of the left graph}
            %\label{}
        \end{subfigure}
        \begin{subfigure}{0.4\textwidth}
            \centering
            \includegraphics[width=1\textwidth]{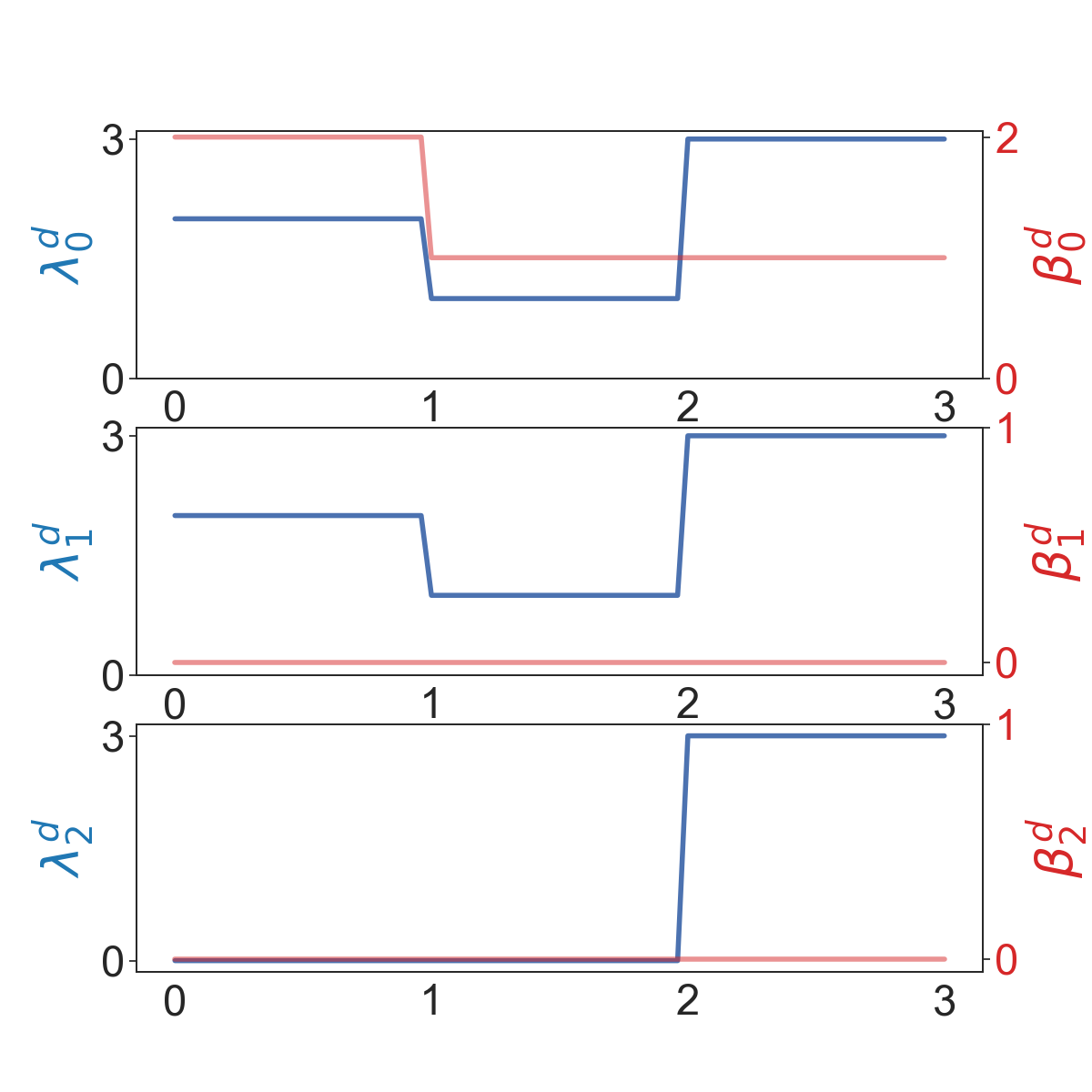}
            \caption{Results of the right graph}
        \end{subfigure}
        \caption{In (a) and (b), numbers on edges are weights. In (c) and (d), the $x$ axis represents the weight. As usual, $\lambda$ and $\beta$ represent the minimal nonzero eigenvalues and Betti numbers.}
   \end{figure}
\end{example}

\subsection{Hypergraphs and hyperdigraphs}

A hypergraph $H$ is a pair $(V, E)$ where $E$ is a subset of the power set of $V$. 
An element $e\in E$ that consists of $q+1$ elements is called a $q$-hyperedge.
To define a chain complex for hypergraphs, the problem is identical to what we encounter in path homology.
If we define the $q$-th chain group to be the vector space generated by $q$-hyperedges, 
the boundary map is not well-defined.
One solution is to consider the \emph{associated simplicial complex} (\emph{simplicial closure}) of 
a hypergraph \cite{parks1991homology}, i.e., the minimal simplicial complex that contains a hypergraph. 
Another solution inspired by the path homology is the \emph{embedded homology} \cite{bressan2019embedded}. 
If we examine the chain complex of the associated simplicial complex, each simplicial chain group $C_q$ contains $D_q$, the vector space generated by $q$-hyperedges. 
We only need to restrict the domain of the simplicial boundary operator to
\begin{align*}
    \text{Inf}_q = D_q \cap \partial_q^{-1}(D_q),
\end{align*}
and then the boundary operator is well-defined.
% the \emph{infimum chain complex}.
\begin{figure}[htbp]
    \centering    
    \begin{tikzcd}
    \dots \arrow[r, "\partial_{q+2}"] & 
    C_{q+1} \arrow[r, "\partial_{q+1}"] & 
    C_{q}  \arrow[r, "\partial_q"] & 
    C_{q-1} \arrow[r, "\partial_{q-1}"] & 
    \dots \\ 
     & 
     D_{q+1} \arrow[u, hook, dashed] &
     D_q \arrow[u, hook, dashed] &
     D_{q-1} \arrow[u, hook, dashed] & 
     & \\  
     & 
     D_{q+1}\cap \partial_{q+1}^{-1}(D_q) \arrow[u, hook, dashed] &
     D_q\cap \partial_{q}^{-1}(D_{q-1}) \arrow[u, hook, dashed] &
     D_{q-1}\cap \partial_{q-1}^{-1}(D_{q-2}) \arrow[u, hook, dashed] & 
     &
    \end{tikzcd}
    \caption{$C_q$ is the $q$-th chain group of the associated simplicial complex of $H$, and $D_q$ is the vector space generated by $q$-hyperedges.}
\end{figure} 

A hyperdigraph is a hypergraph where each hyperedge is ordered\footnote{There are other definitions of a hyperdigraph \cite{ausiello2017directed, thakur2009linear}.},
and its embedded homology can be defined analogously \cite{chen2023persistent}.
Persistent homology of hypergraphs and hyperdigraphs were studied in \cite{bressan2019embedded, ren2018hodge, ren2020stability}.
Persistent hypergraph Laplacians were proposed by Liu et al. \cite{liu2021persistent} and persistent hyperdigraph Laplacians were introduced by Chen et al \cite{chen2023persistent}. 
Alternative approaches regarding the homology and Laplacians of hypergraphs includes \cite{chung1992laplacian, emtander2009betti, hu2015laplacian, jost2019hypergraph, muranov2022path, myers2023topological}.
%Hypergraph has been applied in molecular biology \cite{liu2021persistent, liu2021hypergraph}.

\subsection{Persistent Dirac operators}

Besides Laplacians, Dirac operators on chain complexes have also been studied \cite{ameneyro2024quantum,ameneyro2023quantum, bianconi2021topological, suwayyid2024persistent,suwayyid2024persistentB, wee2023persistent}.
Given a chain complex $(V, d)$
\[
    \begin{tikzcd}[column sep = large]
        \centering
    \cdots \arrow[r, "d_3"] & 
    V_2 \arrow[r, "d_2"] & 
    V_1  \arrow[r, "d_1"] & 
    V_0 \arrow[r] & 
    0
    \end{tikzcd} 
\]
where each chain group $V_q$ is a finite-dimensional inner product space, 
the $q$-th Dirac operator $D_q$ is represented by the block matrix 
\begin{align*}
    \bordermatrix{
        ~       & V_0  & V_1 & V_2 & \cdots & V_q & V_{q+1}\cr
        V_0     & 0   & [d_1] & 0         & \cdots & 0 & 0  \cr
        V_1     & [d_1^{\ast}] & 0 & [d_2] & \cdots & 0 & 0 \cr
        V_2     & 0 & [d_2^{\ast}] & 0 & \cdots & 0 & 0  \cr
        \vdots   & \vdots & \vdots & \vdots & \ddots & \vdots & \vdots \cr
        V_q     & 0 & 0 & 0 & \cdots & 0 & [d_{q+1}] \cr
        V_{q+1} & 0 & 0 & 0 & \cdots & [d_{q+1}^{\ast}] & 0 \cr
    }
\end{align*}where $[]$ denotes a matrix representation of a linear morphism.
Dirac operators are closely related to combinatorial Laplacians.
If we think of all combinatorial Laplacians as a single operator $dd^{\ast} + d^{\ast}d = (d+d^{\ast})^2$ on $V$,
then the $q$-th Dirac operator is the restriction of the square root $d+d^{\ast}$ on $V_0 \oplus \cdots \oplus V_{q+1}$.
We can also see this by direct computation. The square of $D_q$ is 
\begin{align*}
    \bordermatrix{
        ~       & V_0  & V_1 & V_2 & \cdots & V_q & V_{q+1}\cr
        V_0     & [\Delta_0]  & 0 & 0         & \cdots & 0 & 0  \cr
        V_1     & 0 & [\Delta_1] & 0 & \cdots & 0 & 0 \cr
        V_2     & 0 & 0 & [\Delta_2] & \cdots & 0 & 0  \cr
        \vdots   & \vdots & \vdots & \vdots & \ddots & \vdots & \vdots \cr
        V_q     & 0 & 0 & 0 & \cdots & [\Delta_q] & 0 \cr
        V_{q+1} & 0 & 0 & 0 & \cdots & 0 & [\Delta_{q+1,-}] \cr
    }
\end{align*}
where $\Delta_q$ is the $q$-th combinatorial Laplacian. Hence, the square of any eigenvalue $\lambda$ of a Dirac operator must be an eigenvalue of a combinatorial Laplacian.

Recall that when we define persistent Laplacians, we construct an auxiliary subspace $C_{q+1}^{X,Y}$ of $C_{q+1}(Y)$ and a map $\partial_{q+1}^{X,Y}: C_{q+1}^{X,Y}\to C_q(X)$.
Since $C_q(X)$ is actually a subspace of $C_q^{X,Y}$, all $C_q^{X,Y}$ and $\partial_{q}^{X,Y}$ constitute an auxiliary chain complex 
\[
    \begin{tikzcd}[column sep = large]
        \centering
    \cdots \arrow[r, "\partial_3^{X,Y}"] & 
    C_2^{X,Y} \arrow[r, "\partial_2^{X,Y}"] & 
    C_1^{X,Y} \arrow[r, "\partial_1^{X,Y}"] & 
    C_0^{X,Y} \arrow[r] & 
    0
    \end{tikzcd} 
\] The $q$-th persistent Dirac operator of simplicial complexes $X\subset Y$ is just the $q$-th Dirac operator on this auxiliary complex.
The square of a persistent Dirac operator is not necessarily a block matrix consisting of persistent Laplacians. 
It is also possible to extend persistent Dirac operators to other settings such as path complexes and hypergraphs \cite{suwayyid2024persistent}.

\subsection{Mayer homology}

In classical homology theory,  the square of the boundary operator $d$ of a chain complex must be zero ($d^2=0$).  
However, this constraint can be relaxed in the so called Mayer homology theory using an $N$-chain complex \cite{mayer1942new}. 
An $N$-chain complex is a sequence of abelian groups and group morphisms $(V, d)$ where $d^N  = 0$.
A simplicial complex can actually give rise to an $N$-chain complex. 
Recall that in a simplicial chain complex the boundary operator is given by 
\begin{align*}
    \partial [v_{a_0}, \dots, v_{a_q}] = \sum_i (-1)^i[v_{a_0}, \dots, \hat{v}_{a_i}, \dots, v_{a_q}].
\end{align*}
For a prime number $N$, let $\xi = e^{2\pi \sqrt{-1}/N}$, we can define a generalized boundary operator $d$ by 
\begin{align*}
    d [v_{a_0}, \dots, v_{a_q}] = \sum_i \xi^{i} [v_{a_0}, \dots, \hat{v}_{a_i}, \dots, v_{a_q}]
\end{align*}
and prove that $d^N=0$.
Even though $N$-chain complex is not a chain complex in general, observe that for any positive integer $n<N$
\[
    \begin{tikzcd} 
        C_{q+N-n}(X; \cplx) \arrow[r, "d^{N-n}"] & 
        C_{q}(X; \cplx)  \arrow[r, "d^{n}"] & 
        C_{q-n}(X; \cplx) & 
    \end{tikzcd}   
\] resembles a part of chain complex. We can define the Mayer homology group $H_{q, n}(X) = \ker d^{n}/\im d^{N-n}$ \cite{mayer1942new}
and Mayer Laplacians (which can be thought of as $(d^n+(d^{N-n})^{\ast})((d^n)^{\ast}+d^{N-n})$) analogously \cite{shen2023persistent}.
When $N=2$, an $N$-chain complex reduces to a chain complex and 
a Mayer homology group reduces to a normal homology group. 
Shen et al. \cite{shen2023persistent} also introduced persistent Mayer homology and persistent Mayer Laplacians on $N$-chain complexes.  
Suppose $X \subset Y$, then we have the following commutative diagram
\[
    \begin{tikzcd}[column sep = large]
        C_{q+N-n}(X; \cplx) \arrow[r, "(d^X)^{N-n}"] \arrow[d, hook, dashed] 
          & C_{q}(X; \cplx) \arrow[r, "(d^X)^{n}"] \arrow[d, hook, dashed] 
            & C_{q-n}(X; \cplx) \arrow[d, dashed, hook]
            \\
        C_{q+N-n}(Y; \cplx) \arrow[r, "(d^Y)^{N-n}"] 
          & C_q(Y; \cplx) \arrow[r, "(d^Y)^{n}"] 
            & C_{q-n}(Y; \cplx)
    \end{tikzcd}
\] and persistent Mayer Laplacians can be defined analogously.
Compared to simplicial homology,
Mayer homology and Mayer Laplacians provide more features since we can vary parameters $N$ and $n$.
Mayer homology and Mayer Laplacians also concern general relationship between different dimensions. 
Persistent Mayer homology has been applied to protein-ligand binding affinity predictions \cite{feng2024mayer}.

\section {Conclusion and outlook}

The development of persistent topological Laplacians (PTLs) was in a large part inspired by persistent homology's limitations in modeling complex biomolecular data \cite{chen2019evolutionary, wang2020persistent}. 
The techniques we reviewed in this survey transform a point cloud or a network to a set of algebraic features that encode both spatial and non-spatial information. 
The resulting low dimensional features can be further used in supervise or unsupervised machine learning to reveal hidden information. 

% problems, future development
The field of PTLs is dynamic and rapidly evolving. The future development of PTLs is widely open, and we envision the exploration of  the following topics.

(1) Given a point cloud or a network, the diversity of PTLs makes it challenging to fully understand the relationship 
between the geometry and topology of the data and the spectra of PTLs.  
The understanding of this relationship is crucial for the effective application of PTLs to real world problems.  

(2) To a certain extent, the success of TDA can be attributed to its integration with machine learning, 
particularly with the introduction of topological deep learning \cite{cang2017topologynet}. 
Sheaf neural networks \cite{hansen2020sheaf}, sheaf attention \cite{barbero2022sheaf}, and  neural sheaf diffusion \cite{bodnar2022neural} are popular topics. 
Similarly, the development of efficient PTL representations for machine learning, including deep learning, 
is also an important topic.  The featurization of Laplacians typically requires domain knowledge and experience. 
Since self learning representations of persistent diagrams have been proposed \cite{hofer2019learning},
we wonder if self learning representations of (persistent) Laplacians are possible.
As the eigenvectors of PTLs were found to have better descriptive power than eigenvalues \cite{jones2023persistent},
featurization of the eigenvectors of Laplacians is also an interesting future topic.

(3) Despite efforts in computational algorithm and software development \cite{memoli2022persistent, wang2021hermes}, 
the computation of PTLs remains slow,  particularly for problems involving large datasets. 
Since the primary value of TDA lies in its ability to tackle challenges in data science, 
one of the most pressing needs will be the development of efficient and robust PTL software packages to solve real-world problems. 
The development of finite field PTLs will be extremely valuable for the implementation of PTLs in data science.

(4) PTLs have been formulated on a variety of  mathematical  objects, including simplicial complexes, flag complexes,  path complexes, cellular sheaf, digraphs,  hypergraphs, and hyperdigraph. One can also extend PTLs to settings such as the Hochschild complex \cite{gerstenhaber1995higher},  quantum homology \cite{biran2009lagrangian}, multiparameter persistent homology \cite{harrington2019stratifying},  and interaction homotopy and interaction homology \cite{liu2023interaction,liu2024persistent}.  
We expect these developments will further extend the scope and capability of the current TDA for real world applications. 

(5) As mentioned in \cite{suwayyid2024persistent}, a few more persistent Dirac operators can be  formulated for flag complexes, digraphs, hyperdigraphs, etc.
It is possible that persistent sheaf Dirac operators can be devised to distinguish certain geometric shapes. 
Additionally, persistent Dirac operators defined on a spinor bundle may extend persistence to index theory, such as multiscale index theory.  

(6) The PTLs on manifolds, 
such as the evolutionary de Rham-Hodge theory pose implementation challenges compared to their discrete counterparts on point clouds \cite{ribando2024combinatorial}. 
Recently, persistent de Rham-Hodge Laplacian in the Eulerian representation has been proposed for manifold topological learning (MTL) on \cite{su2024persistent}.  
Persistent de Rham-Hodge Laplacians extend earlier persistent homology in cubical setting \cite{kaji2020cubical, wagner2011efficient,wang2016object}. 
From a theoretical point of view, it will be interesting to extend various PTLs to the setting of manifolds (with boundaries). 

(7) In addition to point cloud data and data on manifolds, there are knot-type data, such as DNA packaging in Hi-C data and entangled brain neurons.  
Knots are traditionally studied with invariants, such as the Alexander polynomial, the Jones polynomial, and the Kauffman polynomial \cite{kauffman1987state}. 
Song et al. proposed the multiscale Jones polynomial and the persistent Jones polynomial \cite{song2024multiscale}.
Khovanov homology \cite{khovanov2000categorification} is a major breakthrough in knot theory. Shen et al. \cite{shen2024evolutionary} proposed evolutionary Khovanov homology for weighted links. 
Jones and Wei \cite{jones2024khovanov} proposed Khovanov Laplacians and showed that, at least for chiral prime knots up to 10 crossings, they can distinguish chiral knots from their mirrors. Based on these developments, PTLs on knot-type data, i.e., persistent Khovanov Laplacians, can be formulated, and future research on computational geometric topology is widely open.  

(8) Finally, ChatGPT ushers in a new era of artificial intelligence (AI), offering widespread opportunities in all disciplines. ChatGPT and other chatbots effectively transform pure mathematical theories into practical computational tools, including  PTLs \cite{liu2023chatgpt}. Both AI-enabled topology and topology-enabled AI will have a growing impact on topology\cite{chen2024multiscale}. 
 
\section*{Acknowledgments}
This work was supported in part by NIH grant  R35GM148196, National Science Foundation grants DMS2052983 and IIS-1900473,    MSU Foundation,  and Bristol-Myers Squibb 65109. The authors thank Fei Han,  Fengchun Lei,  Fengling Li, Zhi Lu,  Jie Wu, and Kelin Xia for useful discussions.  

\clearpage
\bibliographystyle{abbrv}
\bibliography{xiaoqi}

\clearpage 
%\section{Appendix}

%\subsection{The calculation of a persistent Laplacian}

\end{document}